\documentclass[11pt]{amsart}
\usepackage{epsfig}
\usepackage{amsmath, amssymb, euscript,  enumerate}
\usepackage{pxfonts}

\newcommand{\dist}{{\text{\rm dist}}}

\newcommand{\ap}{\alpha}             
\newcommand{\bt}{\beta}
\newcommand{\gm}{\gamma}             \newcommand{\Gm}{\Gamma}
\newcommand{\dt}{\delta}             
\newcommand{\vep}{\varepsilon}

\newcommand{\ld}{\lambda}            \newcommand{\Ld}{\Lambda}
\newcommand{\sm}{\sigma}             
\newcommand{\vp}{\varphi}
\newcommand{\om}{\omega}             \newcommand{\Om}{\Omega}
\newcommand{\vr}{\varrho}            \newcommand{\iy}{\infty}
            
\newcommand{\f}{\frac}             \newcommand{\el}{\ell}

\newcommand{\fF}{{\mathfrak F}}

\newcommand{\fL}{{\mathfrak L}}

\newcommand{\fS}{{\mathfrak S}}

\newcommand{\fm}{{\mathfrak m}}

\newcommand{\BN}{{\mathbb N}}

\newcommand{\BR}{{\mathbb R}}

\newcommand{\cC}{{\mathcal C}}

\newcommand{\cI}{{\mathcal I}}
\newcommand{\cJ}{{\mathcal J}}
\newcommand{\cK}{{\mathcal K}}
\newcommand{\cL}{{\mathcal L}}
\newcommand{\cM}{{\mathcal M}}

\newcommand{\Qtil}{{\widetilde Q}}

\newcommand{\s}{\setminus}         \newcommand{\ep}{\epsilon}
\newcommand{\n}{\nabla}            \newcommand{\e}{\eta}
\newcommand{\pa}{\partial}        \newcommand{\fd}{\fallingdotseq}
       
    \newcommand{\ds}{\displaystyle}

 \newcommand{\pf }{\noindent{\it Proof. }}
\newcommand{\rk }{\noindent{\it Remark. }}

\newcommand{\aee }{\text{\rm a.e.}} \newcommand{\diam }{\text{\rm diam}}
  \newcommand{\pv }{\text{\rm p.v.}}
\newcommand{\osc }{\text{\rm osc}}  \newcommand{\dd }{\text{\rm d}}

\newcommand{\rB }{{\text{\rm B}}}   \newcommand{\rC }{{\text{\rm C}}}

\newtheorem{thm}[subsubsection]{Theorem}
\newtheorem{lemma}[subsubsection]{Lemma}
\newtheorem{cor}[subsubsection]{Corollary}
\newtheorem{remark}[subsubsection]{Remark}

\newtheorem{definition}[subsubsection]{Definition}

\newtheorem{assumption}[subsubsection]{Assumption}
\numberwithin{equation}{subsection}

\title[fully nonlinear integro-differential operators]{ Regularity results for fully nonlinear integro-differential operators
with\\ nonsymmetric positive kernels : Subcritical Case }
\author{ Yong-Cheol Kim and Ki-Ahm Lee }
\begin{document}
\begin{abstract}
 We introduce a new class of fully nonlinear
integro-differential operators with possible nonsymmetric kernels.
 For the index $\sigma$ of the
operator  in $ (1,2)$ (subcritical case),  we introduce very
general class of fully nonlinear integro-differential operators and obtain a
comparison principle, a nonlocal version of the
Alexandroff-Backelman-Pucci estimate, a Harnack inequality, a
H\"older regularity, and an interior $\rm C^{1,\alpha}$-regularity
for equations associated with such a class.
\end{abstract}
\address{$\bullet$ Yong-Cheol Kim : Department of Mathematics Education, Korea University, Seoul 136-701,
Korea }

\email{ychkim@korea.ac.kr}

\address{$\bullet$ Ki-Ahm Lee : Department of Mathematics, Seoul National University, Seoul 151-747,
Korea} \email{kiahm@math.snu.ac.kr}

\maketitle

\section{Introduction}\label{sec-1}

It is well-known from general theory on semigroups that the
infinitesimal generator of any L\`evy process always exists for all
functions in the Schwartz space $\fS(\BR^n)$. From the celebrated
L\`evy-Khintchine formula, we can derive that the infinitesimal
generator is given by an operator of the general form
\begin{equation}\label{eq-1.1}
\begin{split}
\cL u(x)&=\sum_{i=1}^n\sum_{j=1}^n a_{ij}\pa_{ij}u+\sum_{i=1}^n
b_i\pa_i u\\&\qquad+\int_{\BR^n\s\{0\}}\bigl[u(x+y)-u(x)-(\n
u(x)\cdot y)\chi_{B_1}(y)\bigr]\,\dd\fm(y).
\end{split}\end{equation}
The first term corresponds to the diffusion, the second to the
drift, and the third to the jump part.

In this paper, we shall focus on the operators which are obtained in
purely jump processes, i.e. processes without diffusion or drift
part. In particular, we are mainly concerned with the operators with
the general form
\begin{equation}\label{eq-1.2}
\cL^t u(x)=\int_{\BR^n\s\{0\}}\bigl[u(x+y)-u(x)-(\n u(x)\cdot
y)\chi_{B_t}(y)\bigr]\,\dd\fm(y),\,\,t>0,
\end{equation}
 where $\fm$ is
a positive measure satisfying $\int_{\BR^n\s\{0\}}(|y|^2\wedge
1)\,\dd\fm(y)<\iy,$ where $|y|^2\wedge 1=\min\{|y|^2,1\}$. The value
of $\cL^t u(x)$ is well-defined when $u$ is bounded on $\BR^n$ and
$\rC^{1,1}$ at $x$. These concepts shall be defined more precisely
later. The operator $\cL^t$ described above is called a {\it linear
integro-differential operator}.

The operator \eqref{eq-1.2} was introduced with too much generality.
So we shall restrict our attention to the operators $\cL^t$ where
the measure $\fm$ is given by a positive kernel $K$ which is not
necessarily symmetric. That is to say, the operators $\cL^t$ are
given by
$$\cL^t u(x)=\pv\int_{\BR^n}\bigl[u(x+y)-u(x)-(\n u(x)\cdot
y)\chi_{B_t}(y)\bigr]K(y)\,dy,\,t>0.$$ Then we see that $\cL^t u(x)$
is well-defined provided that $u\in\rC^{1,1}(x)\cap\rB(\BR^n)$
(refer to Definition \ref{def-2.3} for the definition of
$\rC^{1,1}(x)$) where $\rB(\BR^n)$ denotes {\it the family of all
real-valued bounded functions defined on $\BR^n$}.

If $K$ is symmetric, then an odd function $\bigl[(\n u(x)\cdot
y)\chi_{B_t}(y)\bigr]K(y)$ will be canceled in the integral, and so
we have that
$$\cL^t  u(x)=\pv\int_{\BR^n}\bigl[u(x+y)+u(x-y)-2u(x)\bigr]K(y)\,dy,\,t>0.$$
On the other hand, if $K$ is not symmetric, the effect of $\bigl[(\n
u(x)\cdot y)\chi_{B_t}(y)\bigr]K(y)$ persists and we can actually
observe that the influence of this gradient term becomes stronger as
we try to get an estimate in smaller regions. The aim of this work
is to obtain regularity results for fully nonlinear
integro-differential equations with possible nonsymmetric kernels.

This kind of equations are often obtained in stochastic control
problems \cite{S}. If a player in a stochastic game is permitted to
choose different strategies at every step in order to maximize the
expected value of some function at the first exit point of a domain,
then a convex nonlinear operator
\begin{equation}\label{eq-1.3}
\cI^t u(x)=\sup_{\ap}\cL^t_{\ap}u(x),\,\,t>0,
\end{equation}
appears in the game. In a
competitive game with two or more players, more complicated
operators of the form
\begin{equation}\label{eq-1.4}
\cI^t u(x)=\inf_{\bt}\sup_{\ap}\cL^t_{\ap\bt}u(x),\,\,t>0,
\end{equation}
 can
be obtained. The different aspect between \eqref{eq-1.3} and
\eqref{eq-1.4} is convexity. Also an operator like $\cI^t
u(x)=\sup_{\ap}\inf_{\bt}\cL^t_{\ap\bt}u(x)$ can be considered.
Characteristic properties of these operators can easily be derived
as follows;
\begin{equation}\label{eq-1.5}
\begin{split}\inf_{t\ge 1/2}\inf_{\ap\bt}\cL^t_{\ap\bt}
v(x)&\le\cI^{\pm} [u+v](x)-\cI^{\pm} u(x)\le\sup_{t\ge
1/2}\sup_{\ap\bt}\cL^t_{\ap\bt} v(x),\end{split}
\end{equation}
where $\cI^{\pm} u$ are given by $\cI^+ u(x)=\sup_{t\ge 1/2}\cI^t
u(x)$ and $\cI^- u(x)=\inf_{t\ge 1/2}\cI^t u(x)$ (here we use $1/2$
on $t$ instead of $1$ with certain technical reasons which can be
found in the proof of \ref{thm-7.9} below).

A more general and better description of the fully nonlinear
operators we want to deal with is the operator $\cJ$ for which
\eqref{eq-1.5} holds for some family of linear integro-differential
operators $\cL^t_{\ap\bt}$. The idea is that the concept of
ellipticity can be replaced by an estimate $\cM^-_{\ld,\Ld}
v(x)\le\cJ[u+v](x)-\cJ u(x)\le\cM^+_{\ld,\Ld}v(x)$, where
$\cM^-_{\ld,\Ld}$ and $\cM^+_{\ld,\Ld}$ are the Pucci extremal
operators \cite{CC}. Then it is easy to see that $\cJ$ must be an
elliptic second order differentiable operator. If instead we compare
with suitable nonlocal extremal operators, we will have a concept of
ellipticity for nonlocal operators $\cJ$. We shall give a precise
definition in Section \ref{sec-3} (see Definition \ref{def-3.1}).

We now explain the natural Dirichlet problem for such nonlocal
operators $\cI^{\pm}$. Let $\Om$ be an open domain in $\BR^n$. Given
a function $g$ defined on $\BR^n\s\Om$, we want to find a function
$u$ such that
$$\begin{cases}\cI^{\pm} u(x)=0 &\text{ for any $x\in\Om$,}\\u(x)=g(x) &\text{
for $x\in\BR^n\s\Om$.}\end{cases}$$ Note that the boundary condition
is given not only on $\pa\Om$ but also on the whole complement of
$\Om$. This is because of the nonlocal character of the operators
$\cI^{\pm}$. From the stochastic point of view, it corresponds to
the fact that a discontinuous L\`evy process can exit the domain
$\Om$ for the first time jumping to any point in $\BR^n\s\Om$.

In this paper, we shall concentrate mainly upon the regularity
properties of viscosity solutions to an equation $\cI^{\pm} u(x)=0$.
We shall briefly give a very general comparison principle from which
existence of the solutions can be obtained in smooth domains. Since
kernels of integro-differential operators are comparable to the
kernel of the fractional Laplace operator $-(-\Delta)^{\sm/2}$, the
theory we want to develop can be understood as a theory of viscosity
solutions for fully nonlinear operators of fractional order.

\subsection{The differences between local and nonlocal operators}
Most of all, let us emphasize the main differences between local and nonlocal
equations in terms of Harnack estimate and H\"older continuity.

$\bullet$ First, the nonnegativity of the solution of a local
uniformly elliptic equation in a ball $B_R$ is enough to get a
Harnack inequality in a smaller ball $B_{R/2}$. For the nonlocal
equation, there is a counterexample \cite{BK1} that is nonnegative
in $B_R$ but that has zero value in $B_{R/2}$. It is due to the fact
that there are influence for the values outside $B_R$, which should
be controlled to have a Harnack inequality. So we impose that the
solution is nonnegative in $\BR^n$.

$\bullet$ Second, by applying Harnack inequalities on $\sup_{B_R}
u-u(x)$ and $u(x)-\inf_{B_R}u$, we are able to show the oscillation
lemma and H\"older regularity of the solution of a  local equation.
However such method cannot be directly applicable due to the fact
that we need nonnegativity of the solution in $\BR^n$, not only
$B_R$.

$\bullet$ The last interesting fact is that there is a discontinuous
solution if the Kernel is allowed to oscillate between two different
exponents $\sigma_1$ and $\sigma_2$, \cite{BBC}. The authors show
that there is still a Harnack estimate if the radius of the ball has
a positive lower bound, for example one.

\subsection{History and New results}
There are some known results about Harnack inequalities and H\"older
estimates for integro-differential operators with positive symmetric
kernels (see \cite{J} for analytical proofs and \cite{BBC},
\cite{BK1}, \cite{BK2},\cite{BL}, \cite{KS}, \cite{SV} for
probabilistic proofs). The estimates in all these previous results
blow up as the index $\sm$ of the operator approaches $2$. In this
respect, they do not generalize to elliptic partial differential
equations. However there is some known result on regularity results
for fully nonlinear integro-differential equations associated with
nonlinear integro-differential operators with positive symmetric
kernels which remain uniform as the index $\sm$ of the operator
approaches $2$ (see \cite{CS}). Therefore these results make the
theory of integro-differential operators and elliptic differential
operators become somewhat unified. There has been serious
consideration on the concept of viscosity solution of fully
nonlinear integro-differential equations and their properties,
\cite{ Ar,AT, BCI2,BS,I, JK, P}.

In this paper, we extend the important regularity results of
Caffarelli and Silvestre \cite{CS} on positive symmetric kernels with
certain decay to those on certain positive (not necessarily
symmetric) kernels including such symmetric kernels which remain
uniform as the index $\sm$ of the operator approaches $2$. In this
occasion we can not expect any cancelation on the estimates contrary
to the case of the symmetric kernel, which stirs up some difficulties
in this problem.

Throughout this paper we would like to briefly present the
necessary definitions and then prove some regularity estimates.
Our results in this paper are the following.

$\bullet$ We introduce more general new class of operators to consider
nonsymmetric case. It  is invariant under translation and
scaling, which are crucial properties used at \cite{CS}. Still we
are able to show  standard porperties for viscosity solutions of the general nonlinear
integro-differential equations.

$\bullet$ We show a new version of the nonlocal  Alexandroff-Backelman-Pucci
estimate for fully nonlinear integro-differential equations with
possible nonsymmetric kernel. It contains the extra term caused by the
nonsymmetry of kernel.

$\bullet$ We also give new proof for the construction of special
functions to handle the nonsymmetry of the kernel. And then we show
A Harnack inequality, H\"older regularity , and an interior
$\rC^{1,\ap}$-regularity result for certain fully nonlinear integro-differential equations.

 \subsection{Key Observations}
Nonsymmetric case developes the following differences from the
symmetric case  and nontrivial difficulties.
Key observations are the following:

$\bullet$  For the nonsymmetric case, $K(y)$ and $K(-y)$ can be
chosen any of
  $\ld/|y|^{n+\sm}$ or
  $\Ld/|y|^{n+\sm}$. Therefore there could be an extra
  term $\ds\int_{\BR^n}\frac{\left|(\n u(x)\cdot y)\chi_{B_t}(y)\right|}{|y|^{n+\sigma}}dy$.

$\bullet$ The equation  is not scaling invariant due to
$|\chi_{B_t}(y)|$.

$\bullet$ Somehow the equation has a drift term, not only the diffusion term.
                The case $1<\sm<2$ and the case $0<\sm\le 1$ require different technique
                 due to the difference of the blow rate as $|y|$ approaches to zero and
                 the decay rate as $|y|$ approaches to infinity.
 When $1<\sm<2$, a controllable decay rate of kernel allows H\"older regularity
                     in a larger class, which is invariant under an one-sided
                scaling i.e. if $u$ is a solution of the homogeneous
                equation, then so is  $u_{\ep}(x)=\ep^{-\sm}u(\ep
                x)$ for $0<\ep\le 1$.
                Critical case ($\sm=1$) and supercritical case
                ($0<\sm<1$) have been studied in \cite{KL} with different
                techniques due to the slow decay rate of the kernel as
                $|x|\rightarrow \infty$.
It is noticeable that a gradient term has be considered as a lower
order term ($1<\sigma<2$), \cite{I},  while our gradient effect comes from
the diffusion term which is a main order term.
\subsection{Outline of Paper}
The outline of the paper is as follows. In Section \ref{sec-2},  the
appropriate definitions of subsolutions and supersolutions of fully
nonlinear integro-differential equation in the viscosity sense shall
be given. In the definitions, we shall allow some kind of
discontinuities outside of the domain of the equation. In Section
\ref{sec-3},  we
introduce more general concept of  fully
nonlinear integro-differential equation which is invariant under
one-sided scaling . We define a nonlocal
elliptic operator by comparing its increments with a suitable
maximal operator. This definition is more general than
\eqref{eq-1.4}.  In Section \ref{sec-4}, we study the stability
properties of viscosity solutions given in the definition. A
comparison principle is proven in Section \ref{sec-5} under very
mild assumptions. We guess that one of the most nontrivial results
in the paper is a nonlocal version of the
Alexandroff-Backelman-Pucci estimate to be shown in Section
\ref{sec-6}. It has an extra term caused by the nonsymmetry of the
kernel. It
leads to regularity results for certain fully nonlinear
integro-differential equations. In Section \ref{sec-7}, we construct
a special function, considering nonsymmetry of kernel, and obtain some pointwise estimates which shall
be helpful in proving H\"older estimates in Section \ref{sec-9}. In
Section \ref{sec-8}, we prove a Harnack inequality which plays an
important role in analysis. We obtain the H\"older estimates in
Section \ref{sec-9}. Finally we show an interior
$\rC^{1,\ap}$-estimates in Section \ref{sec-10}.
\section{Viscosity Solutions}\label{sec-2.0}
\subsection{ Definitions}\label{sec-2}

For our purpose, we shall restrict our attention to the operators
$\cL^t$ where the measure $\fm$ is given by a positive kernel $K$
which is not necessarily symmetric. That is to say, the operators
$\cL^t$ are given by
\begin{equation}\label{eq-2.1}
\cL^t u(x)=\pv\int_{\BR^n}\bigl[u(x+y)-u(x)-(\n u(x)\cdot
y)\chi_{B_t}(y)\bigr]K(y)\,dy,\,\,t>0.
\end{equation}
Then we see that $\cL^t u(x)$ is well-defined provided that
$u\in\rC^{1,1}(x)\cap \rB(\BR^n)$ (refer to Definition \ref{def-2.3}
for the definition of $\rC^{1,1}(x)$) where $\rB(\BR^n)$ denotes
{\it the family of all real-valued bounded functions defined on
$\BR^n$}. To simplify the notation, we write
$\mu_t(u,x,y)=u(x+y)-u(x)-(\n u(x)\cdot y) \chi_{B_t}(y),t>0.$ Then
the expression for $\cL^t$ may shortly be written as
\begin{equation}\label{eq-2.2}
\cL^t
u(x)=\int_{\BR^n}\mu_t(u,x,y)K(y)\,dy.
\end{equation}
In particular, for $t>0$, we consider the class $\fL_t$ of the
operators $\cL^t$ associated with the measures $\fm$ given by
positive kernels $K\in\cK_0$ satisfying that
\begin{equation}\label{eq-2.3}(2-\sm)\f{\ld}{|y|^{n+\sm}}\le
K(y)\le(2-\sm)\f{\Ld}{|y|^{n+\sm}},\,\,1<\sm<2.
\end{equation}
In what follows, our main concern shall be on the nonlinear
integro-differential operators which have the form like
\eqref{eq-1.5} where we think that each $\cL_{\ap\bt}^t\in\fL_t$ has
a kernel $K_{\ap\bt}\in\cK_0$ satisfying \eqref{eq-2.3}. The minimum
assumption so that $\cI^{\pm} u$ are well-defined is that every
kernel $K_{\ap\bt}$ must satisfy the following integrability
condition in a uniform way; more precisely, if we set
$K(y)=\sup_{\ap\bt}\,K_{\ap\bt}(y)$, then
\begin{equation}\label{eq-2.4}
\int_{\BR^n}(|y|^2\wedge 1)\,K(y)\,dy<\iy.
\end{equation}
We say that $P$ is a {\it paraboloid of opening $M$} if
\begin{equation}\label{eq-2.5}
P(x)=\el_0+\el(x)\pm\f{M}{2}\,|x|^2
\end{equation}
 where $M$ is a
positive constant, $\el_0$ is real constant and $\el$ is a linear
function. Then $P$ is called {\it convex} when we have $+$ in
\eqref{eq-2.5} and {\it concave} when we have $-$ in \eqref{eq-2.5}.
Let $\Om\subset\BR^n$ be a bounded domain. Given two semicontinuous
functions $u,v$ defined on an open subset $U\subset\Om$ and a point
$x_0\in U$, we say that {\it $v$ touches $u$ by above at $x_0\in U$}
if $u(x_0)=v(x_0)$ and $u(x)\le v(x)$ for any $x\in U$. Similarly,
we say that {\it $v$ touches $u$ by below at $x_0\in U$} if
$u(x_0)=v(x_0)$ and $u(x)\ge v(x)$ for any $x\in U$. For a
semicontinuous function $u$ on $\Om$ and an open subset $U$ of
$\Om$, we define $\Theta^+ (u,U)(x_0)$ to be the infimum of all
positive constants $M$ for which there is a convex paraboloid of
opening $M$ that touches $u$ by above at $x_0\in U$. Also we define
$\Theta^+ (u,U)(x_0)=\iy$ if no such constant $M$ exists. Similarly,
we define $\Theta^- (u,U)(x_0)$ to be the infimum of all positive
constants $M$ for which there is a concave paraboloid of opening $M$
that touches $u$ by below at $x_0\in U$, and also we define
$\Theta^- (u,U)(x_0)=\iy$ if no such constant $M$ exists. Finally we
set $\Theta(u,U)(x_0)=\max\{\Theta^+ (u,U)(x_0),\Theta^-
(u,U)(x_0)\}\le\iy.$ For these definitions, the readers can refer to
\cite{CC}.
\begin{definition}\label{def-2.1}
Let $\Om\subset\BR^n$ be a bounded domain and let $u:\Om\to\BR$ be a
semicontinuous function. Given $x_0\in\Om$, we say that $u$ is
$\rC_{\pm}^{1,1}$ {\rm at $x_0$} $($ resp. $\rC^{1,1}$ {\rm at
$x_0$} $)$ if $\,\Theta^{\pm} (u,U)(x_0)<\iy$ $($ resp.
$\Theta(u,U)(x_0)<\iy$ $)$ for some open neighborhood $U$ of $x_0$
and we write $u\in\rC^{1,1}[x_0]$ if $\,\Theta(u,U)(x_0)<\iy$. Given
a fixed $\ep\in (0,1)$, we set $\Theta(u,\ep)(x)=\Theta(u,\Om\cap
B_{\ep}(x))$ for $x\in\Om$ and we write $u\in\rC^{1,1}[\Om]$ if
$\sup_{x\in\Om}\Theta(u,\ep)(x)\fd\Theta_{\Om}[u]<\iy$.
\end{definition}
\begin{remark}\label{rem-2.2}
\begin{enumerate}[$(a)$]
\item We note that if $u\in\rC^{1,1}[x_0]$ then $u$ is
differentiable at $x_0$ because $u$ lies between two tangent
paraboloids in an open neighborhood of $x_0$.
\item  In fact, the number $\Theta_{\Om}[u]$ in Definition \ref{def-2.1}  depends
upon $\ep$. But existence of $\ep$ so that $\Theta_{\Om}[u]<\iy$ is
enough, because $\Theta_{\Om}[u]$ decreases as $\ep$ approaches $0$.
\end{enumerate}
\end{remark}
\begin{definition}\label{def-2.3}
A function $u:\BR^n\to\BR$ is said to be $\text{\rm $\rC^{1,1}$ at a
point $x\in\BR^n$}$ $($we write $u\in\rC^{1,1}(x)$$)$, if there
exist some vector $v\in\BR^n$, $r_0>0$ and $M>0$ such that
\begin{equation}\label{eq-2.6}
\bigl|u(x+y)-u(x)-v\cdot y\bigr|\le M\,|y|^2\,\,\,\text{ for
any $y\in B_{r_0}$.}
\end{equation}
We write $u\in\rC^{1,1}(U)$ if
$\,u\in\rC^{1,1}(x)$ for any $x\in U$ and the constant $M$ in
\eqref{eq-2.6} is independent of $\,x$, where $U$ is an open subset of
$\BR^n$.
\end{definition}
\begin{remark}\label{rem-2.4}
\begin{enumerate}[$(a)$]
\item  Such vector $v$ exists uniquely and moreover $v=\n u(x)$.
\item  If $\,u\in\rC^{1,1}[x]$ for $x\in\Om$, then we easily see that
$u\in\rC^{1,1}(x)$. Moreover, it is easy to show that the converse
holds. Thus we conclude that $\rC^{1,1}[\Om]=\rC^{1,1}(\Om)$.
\end{enumerate}
\end{remark}
For $x\in\Om$ and a function $u:\BR^n\to\BR$ which is semicontinuous
on $\overline\Om$, we say that $\vp$ belongs to the function class
$\rC^2_{\Om}(u;x)^+$ (resp. $\rC^2_{\Om}(u;x)^-$) and we write
$\vp\in\rC^2_{\Om}(u;x)^+$ (resp. $\vp\in\rC^2_{\Om}(u;x)^-$) if
there are an open neighborhood $U\subset\Om$ of $x$ and
$\vp\in\rC^2(U)$ such that $\vp(x)=u(x)$ and $\vp>u$ (resp. $\vp<u$)
on $U\s\{x\}$. We note that geometrically $u-\vp$ having a local
maximum at $x$ in $\Om$ is equivalent to $\vp\in\rC^2_{\Om}(u;x)^+$
and $u-\vp$ having a local minimum at $x$ in $\Om$ is equivalent to
$\vp\in\rC^2_{\Om}(u;x)^-$. For $x\in\Om$ and $\vp\in
\rC^2_{\Om}(u;x)^{\pm}$, we write
$$\mu_t(u,x,y;\n\vp)=u(x+y)-u(x)-(\n\vp(x)\cdot y)\chi_{B_t}(y),\,\,t>0,$$
and the expression for $\cL_{\ap\bt}^t\,u(x;\n\vp)$ and $\cI^t
u(x;\n\vp)$ may be written as
\begin{equation*}\begin{split}\cL_{\ap\bt}^t\,u(x;\n\vp)&=\int_{\BR^n}\mu_t(u,x,y;\n\vp)K_{\ap\bt}(y)\,dy,\\
\cI^t
u(x;\n\vp)&=\inf_{\bt}\sup_{\ap}\cL_{\ap\bt}^t\,u(x;\n\vp),\end{split}\end{equation*}
where $K_{\ap\bt}\in\cK_0$. We set ${\cI^+}u(x;\n\vp)=\sup_{t\ge
1/2}\cI^t u(x;\n\vp)$ and ${\cI^-}u(x;\n\vp)=\inf_{t\ge 1/2}\cI^t
u(x;\n\vp)$. We note that if $u\in\rC^{1,1}(x)$, then
$\mu_t(u,x,y;\n\vp)=\mu_t(u,x,y)$, $\cL_{\ap\bt}^t
u(x;\n\vp)=\cL_{\ap\bt}^t u(x)$ and $\cI^t u(x;\n\vp)=\cI^t u(x)$
for $x\in\Om$ and $\vp\in\rC^2_{\Om}(u;x)^{\pm}$.
\begin{definition}\label{def-2.5}
Let $f:\BR^n\to\BR$ be a function. Then a function $u:\BR^n\to\BR$
which is upper $($lower$)$ semicontinuous on $\overline\Om$ is said
to be a {\rm viscosity subsolution (viscosity supersolution)} of an
equation $\cJ u=f$ on $\Om$ and we write $\cJ u\ge f$ $($$\cJ u\le
f$$)$ on $\Om$ in the viscosity sense, if for any $x\in\Om$ there is
some open neighborhood $U$ of $x$ with $U\subset\Om$ such that $\cJ
v(x)$ is well-defined and $\cJ v(x)\ge f(x)$ $($$\cJ v(x)\le
f(x)$$)$ for $v=\vp\chi_{U}+u\chi_{\BR^n\s U}$ whenever
$\vp\in\rC^2(U)$ satisfying $\vp(x)=u(x)$ and $\vp>u$ $($$\vp<u$$)$
on $U\s\{x\}$ exists. Also a function $u$ which is both a viscosity
subsolution and a viscosity supersolution of the equation $\cJ u=f$
on $\Om$ is said to be a $\text{\rm viscosity solution}$ to $\cJ
u=f$ on $\Om$.
\end{definition}

\begin{remark}\label{rem-2.6} Definition \ref{def-2.5} is essentially the same
as Definition 2 in \cite{BI}.\end{remark}

Instead of test functions $\vp\in\rC^2$ given in the above,
functions $\vp$ which are $\rC^{1,1}$ only at the contact point $x$
could be used. This is a larger set of test functions, so that a
priori it may provide a stronger concept of solution. In Section \ref{sec-4},
we shall show that the two approaches are actually equivalent.
\begin{thm}\label{thm-2.7}
Let $f:\BR^n\to\BR$ be a function. Then we have the followings:

$(a)$ If $u:\BR^n\to\BR$ is a function which is upper semicontinuous
on $\overline\Om$, then $\cI^{\pm}u\ge f$ on $\Om$ in the viscosity
sense if and only if $\,{\cI^{\pm}}u(x;\n\vp)$ is well-defined and
\begin{equation}\label{eq-2.7}{\cI^{\pm}}u(x;\n\vp)\ge
f(x)\,\,\text{ for any $x\in\Om$ and $\vp\in\rC^2_{\Om}(u;x)^+$.}
\end{equation}

$(b)$ If $u:\BR^n\to\BR$ is a function which is lower semicontinuous
on $\overline\Om$, then $\cI^{\pm}u\le f$ on $\Om$ in the viscosity
sense if and only if $\,{\cI^{\pm}}u(x;\n\vp)$ is well-defined and
\begin{equation}\label{eq-2.8}
{\cI^{\pm}}u(x;\n\vp)\le f(x)\,\,\text{ for any $x\in\Om$ and
$\vp\in\rC^2_{\Om}(u;x)^-$.}
\end{equation}

$(c)$ If $u:\BR^n\to\BR$ is a function which is continuous on
$\overline\Om$, then $u$ is a viscosity solution to $\cI^{\pm} u=f$
on $\Om$ if and only if it satisfies both \eqref{eq-2.7} and
\eqref{eq-2.8}.
\end{thm}
\pf Since (b) and (c) can be obtained by the similar way to the
proof of (a), we have only to prove the equivalence (a) for $\cI^-$;
similarly for $\cI^+$.

Assume that $\cI^- u\ge f$ on $\Om$ in the viscosity sense. Fix any
$x\in\Om$ and take any $\vp\in\rC^2_{\Om}(u;x)^+$. Then there is an
open neighborhood $U\subset\Om$ of $x$ such that $u(x)=\vp(x)$ and
$\vp>u$ on $U\s\{x\}$. For $0<s<\dd(x,\pa U)$, we set
$u_s=\vp\chi_{B_s(x)}+u\chi_{\BR^n\s B_s(x)}$ and
$v=\vp\chi_U+u\chi_{\BR^n\s U}$. Then, for any $t\ge 1/2$ and any
$\bt$, there is some $\ap$ such that $\cL^t_{\ap\bt}
v(x)\ge\cL^t_{\ap\bt}u_s(x)\ge f(x).$ By the Lebesgue's dominated
convergence theorem, taking $s\downarrow 0$ we conclude that for any
$t\ge 1/2$ and any $\bt$ there is some $\ap$ such that
$\cL^t_{\ap\bt}v(x)\ge\cL^t_{\ap\bt}u(x;\n\vp)\ge f(x)$. Therefore
$\cI^- u(x;\n\vp)$ is well-defined and $\cI^- u(x;\n\vp)\ge f(x)$.

Conversely, suppose that \eqref{eq-2.7} holds. Let $u:\BR^n\to\BR$
be a function which is upper semicontinuous on $\overline\Om$. Fix
any $x\in\Om$ and take a function $\vp\in\rC^2(U)$ satisfying that
$\vp(x)=u(x)$ and $\vp>u$ on $U\s\{x\}$ where $U$ is an open
neighborhood $U\subset\Om$ of $x$. Then $\vp\in\rC^2_{\Om}(u;x)^+$.
If we set $v=\vp\chi_U+u\chi_{\BR^n\s U}$, then for any $t\ge 1/2$
and any $\bt$ there is some $\ap$ such that
$\cL^t_{\ap\bt}v(x)\ge\cL^t_{\ap\bt}u(x;\n\vp)\ge f(x)$, and thus
$\cI^- v(x)\ge f(x)$. Hence we conclude that $\cI^- u\ge f$ on $\Om$
in the viscosity sense. \qed

\subsection{ Maximal operators}\label{sec-3}

In \eqref{eq-1.3} and \eqref{eq-1.4}, we considered the supremum or
an inf-sup of a collection of linear integro-differential operators.
Let us consider a class $\fL$ of linear operators which includes the
class $\fL_0=\bigcup_{t\ge 1/2}\fL_t$ given in Section \ref{sec-2}.
The maximal operator and the minimal operator with respect to $\fL$
are defined by
\begin{equation}\label{eq-3.1}\cM^+_{\fL}u(x)=\sup_{\cL\in\fL}\cL u(x)\,\,\text{ and }\,\,\cM^-_{\fL}u(x)=\inf_{\cL\in\fL}\cL u(x).
\end{equation}
For a function $u:\BR^n\to\BR$ semicontinuous on $\overline\Om$ and
$\vp\in\rC^2_{\Om}(u;x)^{\pm}$, we set
$$\cM^+_{\fL_0}u(x;\n\vp)=\sup_{t\ge 1/2}\cM^+_{\fL_t}u(x;\n\vp)\,\text{ and }\,\cM^-_{\fL_0}u(x;\n\vp)=\inf_{t\ge 1/2}\cM^-_{\fL_t}
u(x;\n\vp),$$ where $\cM^+_{\fL_t}u(x)=\sup_{\cL\in\fL_t}\cL u(x)$
and $\cM^-_{\fL_t}u(x)=\inf_{\cL\in\fL_t}\cL u(x).$ Then we see that
$\cM^-_{\fL_0}u(x;\n\vp)\le{\cI^{\pm}}u(x;\n\vp)\le\cM^+_{\fL_0}u(x;\n\vp)$,
and $\cM^{\pm}_{\fL_0}u(x;\n\vp)$ have the following simple forms;
\begin{equation}\label{eq-3.2}
\begin{split}
&\cM^+_{\fL_0}u(x;\n\vp)=(2-\sm)\int_{\BR^n}\f{\Ld\overline{\mu}^+(u,x,y;\n\vp)
-\ld\underline{\mu}^-(u,x,y;\n\vp)}{|y|^{n+\sm}}\,dy,\\
&\cM^-_{\fL_0}u(x;\n\vp)=(2-\sm)\int_{\BR^n}\f{\ld\underline{\mu}^+(u,x,y;\n\vp)
-\Ld\overline{\mu}^-(u,x,y;\n\vp)}{|y|^{n+\sm}}\,dy,
\end{split}
\end{equation}
where $\overline{\mu}^+$, $\underline{\mu}^+$, $\overline{\mu}^-$
and $\underline{\mu}^-$ are given by
\begin{equation*}\begin{split}\overline{\mu}^{\pm}(u,x,y;\n\vp)=\sup_{t\ge
1/2}\mu^{\pm}_t(u,x,y;\n\vp)\text{ and
}\underline{\mu}^{\pm}(u,x,y;\n\vp)=\inf_{t\ge
1/2}\mu^{\pm}_t(u,x,y;\n\vp).\end{split}\end{equation*} We note that
if $u\in\rC^{1,1}(x)$, then
$\cM^{\pm}_{\fL_0}u(x;\n\vp)=\cM^{\pm}_{\fL_0}u(x)$. We shall use
these maximal and minimal operators to obtain regularity estimates.
The factor $(2-\sm)$ is important when $\sm\to 2$, because we need
such factor if we want to obtain second order partial differential
equations as limits of linear integro-differential equations. In
terms of the regularity, we need the factor $(2-\sm)$ for the
estimates not to blow up as $\sm\to 2$.

Let $K(x)=\sup_{\ap}K_{\ap}(x)$ where $K_{\ap}$'s are all the
kernels of all operators in $\fL$. Instead of \eqref{eq-2.4}, for
any class $\fL$ we shall assume that
\begin{equation}\label{eq-3.3}
\int_{\BR^n}(|y|^2\wedge 1)\,K(y)\,dy<\iy.
\end{equation}
 Using
the extremal operators, we provide a general definition of
ellipticity for nonlocal equations. The following is a kind of
operators of which the regularity result shall be obtained in this
paper.

\begin{definition}\label{def-3.1} Let $\fL$ be a class of linear
integro-differential operators. Assume that \eqref{eq-3.3} holds for $\fL$.
Then we say that an operator $\cJ$ is $\text{\rm elliptic with
respect to $\fL$}$, if it satisfies the following properties:

$(a)$ $\cJ u(x)$ is well-defined for any $u\in\rC^{1,1}[x]\cap
\rB(\BR^n)$.

$(b)$ $\cJ u$ is continuous on an open set $\Om\subset\BR^n$,
whenever $u\in\rC^{1,1}[\Om]\cap\rB(\BR^n)$.

$(c)$ If $\,u,v\in\rC^{1,1}[x]\cap\rB(\BR^n)$, then we have that
\begin{equation}\label{eq-3.4}
\cM^-_{\fL}[u-v](x)\le\cJ u(x)-\cJ
v(x)\le\cM^+_{\fL}[u-v](x).
\end{equation}
\end{definition}

We shall show that any operator as in \eqref{eq-1.5} is elliptic
with respect to any class containing all the operators
$\cL^t_{\ap\bt}$ as long as the condition \eqref{eq-2.4} is
satisfied (see Lemma \ref{lem-3.2} and Lemma \ref{lem-4.2}).

\begin{lemma}\label{lem-3.2}
Let $\cI^{\pm}$ be the operators as in \eqref{eq-1.5} so that
\eqref{eq-2.4} holds for every $K_{\ap\bt}$ and let $\fL$ be any
collection of linear integro-differential operators. If
$\,\fL_0\subset\fL$ where $\fL_0=\bigcup_{t\ge 1/2}\fL_t$, then we
have that
\begin{equation}
\begin{split}
\cM^-_{\fL}[u-v](x)&\le\cI^{\pm} u(x)-\cI^{\pm}
v(x)\le\cM^+_{\fL}[u-v](x)
\end{split}
\end{equation}
for any $u,v\in \rC^{1,1}[x]\cap\rB(\BR^n)$.
\end{lemma}

\pf It can be shown in a similar way as in \cite{CS}. \qed

Definition \ref{def-2.3} is not set up to evaluate the operator
$\cJ$ in the original function $u$. Whenever a smooth function $\vp$
touches $u$ from above, we can always construct a test function
$v\in\rC^{1,1}[x]$ to evaluate $\cJ$. It is very interesting that if
$\cJ$ is any nonlinear operator which is formulated by an inf-sup
(or a sup-inf) of linear operators satisfying \eqref{eq-2.3} and
\eqref{eq-2.4}, then without constructing such a test function $\cJ$
can be evaluated classically in $u$ at those points $x$ where $u$
can be touched by above with a $\rC^2$ function. This interesting
fact is proved in the next lemma.

\begin{lemma}\label{lem-3.3} Let $\cI^{\pm}$ be the operator as in
\eqref{eq-1.5} so that \eqref{eq-2.4} holds for every $K_{\ap\bt}$ and let $\,u\in
\rB(\BR^n)$ be a viscosity subsolution to $\cI^{\pm} u=f$ on $\Om$.
If $\,u\in\rC^{1,1}[x]$ for a point $x\in\Om$, then $\cI^{\pm} u(x)$
is defined in the classical sense and $\cI^{\pm} u(x)\ge
f(x)$.\end{lemma}

\pf It can be obtained in a similar way as in \cite{CS}. \qed

In the next theorem, we shall obtain a result on the operators
$\cM^{\pm}_{\fL_0}$ which is similar to Theorem \ref{thm-2.7}. The
proof is almost the same as that of Theorem \ref{thm-2.7}, and so we
will just write out the statement without detailed proof.

\begin{thm}\label{thm-3.4} Let $f:\BR^n\to\BR$ be a function. Then we
have the followings:

$(a)$ If $u:\BR^n\to\BR$ is a function which is upper semicontinuous
on $\overline\Om$, then $u$ is a viscosity subsolution to
$\cM^+_{\fL_0} u=f$ on $\Om$ if and only if $\,\cM^+_{\fL_0}
u(x;\n\vp)$ is well-defined and $\cM^+_{\fL_0} u(x;\n\vp)\ge
f(x)\,\,\text{ for any $x\in\Om$ and any
$\vp\in\rC^2_{\Om}(u;x)^+$.}$

$(b)$ If $u:\BR^n\to\BR$ is a function which is lower semicontinuous
on $\overline\Om$, then $u$ is a viscosity supersolution to
$\cM^-_{\fL_0} u=f$ on $\Om$ if and only if $\,\cM^-_{\fL_0}
u(x;\n\vp)$ is well-defined and $\cM^-_{\fL_0} u(x;\n\vp)\le
f(x)\,\,\text{ for any $x\in\Om$ and any
$\vp\in\rC^2_{\Om}(u;x)^-$.}$
\end{thm}

\subsection{ Stability properties}\label{sec-4}

In this section, we obtain a few technical properties of the
operators $\cI^{\pm}$ as in \eqref{eq-1.5}. First we shall show that
if $u\in \rC^{1,1}[\Om]\cap\rB(\BR^n)$, then $\cI^{\pm} u$ are
continuous on $\Om$. As mentioned in the previous sections, it is
necessary to justify that the operators of the form \eqref{eq-1.3}
and \eqref{eq-1.4} satisfy the conditions of Definition
\ref{def-3.1}. Next we shall show that our notion of viscosity
solutions allows to touch with solutions which are only punctually
$\rC^{1,1}$ instead of $\rC^2$ in a neighborhood of the point. Then
we shall show the important stability property of viscosity
solutions given in Definition \ref{def-2.3}. We now start with
several technical lemmas.

\begin{lemma}\cite{CS}\label{lem-4.1} Let $\{h_{\ap}\}$ be a family of functions
such that $|h_{\ap}(x)|\le h(x)$ for $h\in L^1(\BR^n)$. If $f\in
L^{\iy}(\BR^n)$, then the family $\{f*h_{\ap}\}$ is uniformly
equicontinuous on every compact subsets of $\BR^n$.\end{lemma}

When we gave the definition of viscosity solutions in Section
\ref{sec-2}, we used $\rC^2$ test functions. Now we show that it is
equivalent to use punctually $\rC^{1,1}$ functions.

\begin{lemma}\label{lem-4.2} Let $\{\cI^t\}_{t\ge 1/2}$ be the family of
operators as in \eqref{eq-1.3} and \eqref{eq-1.4} satisfying
\eqref{eq-2.4} and let $\Omega\subset\BR^n$ be a bounded domain. If
$u\in\rC^{1,1}[\Om]\cap \rB(\BR^n)$, then $\cI^{\pm} u$ are
continuous on $\Om$.\end{lemma}

\pf Fix any $t>0$. Then we have only to prove that the family
$\{\cL^t_{\ap\bt} u\}$ is equicontinuous on $\Om$. We set
$K(x)=\sup_{\ap\bt}K_{\ap\bt}(x)$ as in \eqref{eq-2.4}. Take any
$\ep>0$ and $x_0\in\Om$. Since $\Theta_{\Om}[u]<\iy$, we have that
$\bigl|\mu_t(u,x,y)\bigr|<\Theta_{\Om}[u]\,|y|^2$, whenever
$x\in\Om$ and $|y|<t_1\fd\dist(x,\pa\Om)\wedge t$. Then choose some
sufficiently small $t_0\in(0,t_1)$ so that
\begin{equation}\label{eq-4.1}\int_{B_{t_0}}\Theta_{\Om}[u]\,|y|^2 K(y)\,dy<\ep/3.
\end{equation}
Now we have that
\begin{equation}\label{eq-4.2}
\begin{split}\cL^t_{\ap\bt}
u(x)&=\int_{B_{t_0}}\mu_t(u,x,y)K_{\ap\bt}(y)\,dy+\int_{\BR^n\s
B_{t_0}}\mu_t(u,x,y)K_{\ap\bt}(y)\,dy\\&\fd\cL^{t,0}_{\ap\bt}u(x)+\cL^{t,1}_{\ap\bt}u(x).
\end{split}\end{equation}
From \eqref{eq-4.1}, we easily obtain that $|\cL^{t,0}_{\ap\bt}
u(x)|\le\int_{B_{t_0}}\Theta_{\Om}[u]\,|y|^2 K(y)\,dy<\ep/3$ for any
$\ap,\bt$, whenever $x\in\Om$. We also write
\begin{equation}\label{eq-4.3}
\begin{split}\cL^{t,1}_{\ap\bt}u(x)&=
u*h_{\ap\bt}(x)-\biggl(\int_{\BR^n}h_{\ap\bt}(y)\,dy\biggr)u(x)-\int_{B_t}[\n
u(x)\cdot
y]\,h_{\ap\bt}(y)\,dy\\&\fd\cL^{t,2}_{\ap\bt}u(x)+\cL^{t,3}_{\ap\bt}u(x)+\cL^{t,4}_{\ap\bt}u(x)
\end{split}\end{equation}
where $h_{\ap\bt}(y)=K_{\ap\bt}(y)\chi_{\BR^n\s B_{t_0}}(y)$. Since
$\Theta_{\Om}[u]<\iy$ and $\n u$ is Lipschitz continuous on $\Om$ by
standard analysis, we have that
\begin{equation}\label{eq-4.4}\begin{split}\bigl|\cL^{t,4}_{\ap\bt}u(x)-\cL^{t,4}_{\ap\bt}u(x_0)\bigr|
&\le\biggl(\int_{\BR^n}h_{\ap\bt}(y)|y|\,dy\biggr)\,\Theta_{\Om}[u]\,\,|x-x_0|
\end{split}
\end{equation}
whenever $x\in B_s(x_0)$ and $B_s(x_0)\subset\Om$. Since
$\sup_{\ap,\bt}\int_{\BR^n}h_{\ap\bt}(y)|y|\,dy<\iy$ ( by
\eqref{eq-2.4} ) and $u\in\rC^{1,1}[\Om]$, by Lemma \ref{lem-4.1},
\eqref{eq-4.3} and \eqref{eq-4.4} there exists some sufficiently
small $\dt=\dt(\vep)>0$ such that
\begin{equation}\label{eq-4.5}\bigl|\cL^{t,1}_{\ap\bt}u(x)-\cL^{t,1}_{\ap\bt}u(x_0)\bigr|<\ep/3\end{equation}
for any $\ap,\bt$, and $t\ge 1/2$, whenever $x\in\Om$ and
$|x-x_0|<\dt$. Thus it follows from \eqref{eq-4.2} and
\eqref{eq-4.5} that
$$\bigl|\cL^t_{\ap\bt}u(x)-\cL^t_{\ap\bt}u(x_0)\bigr|\le\bigl|\cL^{t,0}_{\ap\bt}u(x)\bigr|
+\bigl|\cL^{t,0}_{\ap\bt}u(x_0)\bigr|+\bigl|\cL^{t,1}_{\ap\bt}u(x)-\cL^{t,1}_{\ap\bt}u(x_0)\bigr|<\ep$$
for any $\ap,\bt$, and $t\ge 1/2$, whenever $x\in\Om$ and
$|x-x_0|<\dt$. Hence this implies that $|\cI^{\pm} u(x)-\cI^{\pm}
u(x_0)|<\ep$, whenever $x\in\Om$ and $|x-x_0|<\dt$. Therefore we
complete the proof. \qed

When we gave the definition of viscosity solutions in Section
\ref{sec-2}, we used $\rC^2$ test functions. We show it is
equivalent to use punctually $\rC^{1,1}$ functions.

\begin{lemma}\label{lem-4.3} Let $\cJ$ be elliptic with respect to a
class $\fL$ in the sense of Definition \ref{def-3.1}. Assume that
$u:\BR^n\to\BR$ is a viscosity subsolution to $\cJ u=f$ on $\Om$ and
$\vp\in\rC^{1,1}[x]\cap\rB(\BR^n)$ for $x\in\Om$. If $\,\vp$ touches
globally $u$ from above at $x$, then $\cJ\vp(x)$ is defined in the
classical sense and $\cJ\vp(x)\ge f(x)$.\end{lemma}

\pf It can be shown in a similar way as in \cite{CS}. \qed

One of the most useful properties of viscosity solutions is their
stability property under uniform limits on compact sets. We shall
prove a slightly stronger result that the notion of viscosity
subsolution (supersolution) is stable with respect to the natural
limits for upper (lower) semicontinuous functions. This type of
limit is well-known and usually called {\it $\Gm$-limit}. It was
originally called as the celebrated half relaxed limit techniques by
Barles and Perthame, but we are going to follow the similar
definitions at \cite{CS}.

\begin{definition}\label{def-4.4} A sequence $\{u_k\}$ of${\,}$ lower
semicontinuous functions is said to $\text{\rm $\Gm$-converges to
$u$}$ on a set $\Om\subset\BR^n$ if
$(a)$ for any sequence $\{x_k\}\subset\Om$ with
$\lim_{k\to\iy}x_k=x$, $\liminf_{k\to\iy}u_k(x_k)\ge u(x)$ and $(b)$
for any $x\in\Om$, there is a sequence $\{x_k\}\subset\Om$ with
$\lim_{k\to\iy}x_k=x$ such that
$\lim_{k\to\iy}u_k(x_k)=u(x)$.\end{definition}

\noindent{\it Remark.} (a) A uniformly convergent sequence $\{u_k\}$
converges in the $\Gm$ sense.

(b) If $\{u_k\}$ $\Gm$-converges to $u$ on $\Om$ and $u$ has a
strict local minimum at $x$ then there is a sequence $\{x_k\}$ with
$\lim_{k\to\iy}x_k=x$ such that each $u_k$ has a local minimum at
$x_k$ (see \cite{GD}).

(c) If $\{u_k\}$ $\Gm$-converges to $u$ on $\Om$, then $\{u_k-\vp\}$
$\Gm$-converges to $u-\vp$ on $\Om$ where $\vp\in\rC(\Om)$.

(d) From (b) and (c), we can get that if $\{u_k\}$ $\Gm$-converges
to $u$ on $\Om$ and $u-\vp$ has a strict local minimum at $x$ where
$\vp\in\rC(\Om)$, then there is a sequence $\{x_k\}$ with
$\lim_{k\to\iy}x_k=x$ such that each $u_k-\vp$ has a local minimum
at $x_k$.

\begin{lemma}\label{lem-4.5} Let $\cJ$ be elliptic in the sense of
Definition \ref{def-3.1}. If $\,\{u_k\}\subset\rB(\BR^n)$ is a
sequence of viscosity supersolutions $u_k$ to $\cJ u_k=f_k$ on $\Om$
such that

$(a)$ $\{u_k\}$ $\Gm$-converges to $u$ in $\Om$, $(b)$ $\{u_k\}$
converges to $u$ $\aee$ on $\BR^n$ and

$(c)$ $\{f_k\}$ converges to $f$ locally uniformly on $\Om$, then
$u$ is a viscosity supersolution to $\cJ u=f$ on $\Om$.\end{lemma}

\pf It can be done with minor changes in a similar way as in
\cite{CS}. \qed

We just obtained the stability property of supersolutions under
$\Gm$-limits. For the corresponding result for subsolutions as in
the following lemma, we would also consider the natural limit in the
space of upper semicontinuous functions which is the same as the
$\Gm$-convergence of $-u_k$ to $-u$.

\begin{lemma}\label{lem-4.6} Let $\cJ$ be elliptic in the sense of
Definition \ref{def-3.1}. If $\,\{u_k\}\subset\rB(\BR^n)$ is a sequence of
viscosity subsolutions to $\cJ u_k=f_k$ on $\Om$ such that

$(a)$ $\{-u_k\}$ $\Gm$-converges to $-u$ in $\Om$, $(b)$ $\{u_k\}$
converges to $u$ $\aee$ on $\BR^n$ and

$(c)$ $\{f_k\}$ converges to $f$ locally uniformly on $\Om$, then
$u$ is a viscosity subsolution to $\cJ u=f$ on $\Om$.\end{lemma}

As a corollary, we also obtain the stability property under
uniform limits.

\begin{cor}\label{cor-4.7} Let $\cJ$ be elliptic in the sense of
Definition \ref{def-3.1}. If $\,\{u_k\}\subset\rB(\BR^n)$ is a sequence of
viscosity solutions to $\cJ u_k=f_k$ on $\Om$ such that

\noindent$\,\,\,\,\,(a)$ $\{u_k\}$ and $\{f_k\}$ converge to $u$ and
$f$ locally uniformly on $\Om$, respectively,

\noindent$\,\,\,\,\,(b)$ $\{u_k\}$ converges to $u$ $\aee$ on
$\BR^n$, then $u$ is a viscosity solution to $\cJ u=f$ on
$\Om$.\end{cor}

\pf Since $u_k\to u$ locally uniformly on $\Om$, we see that
$\{u_k\}$ $\Gm$-converges to $u$ in $\Om$. Thus the required result
follows from Lemma \ref{lem-4.5} and \ref{lem-4.6}. \qed
\subsection{ Comparison principle}\label{sec-5}

The comparison principle for viscosity solutions can be shown by
very standard ideas in nonlinear analysis, which originated from the
idea of Jensen \cite{J} using sup-convolutions and inf-convolutions.
The method has been succesfully adapted to integro-differential
equations \cite{A} and a more general proof can be found in
\cite{BI} in case that the viscosity solutions have an arbitrary
growth at infinity. Our definitions do not quite fit in with the
previous frameworks because we consider mainly the general class of
operators given by Definition \ref{def-3.1} and we allow
discontinuities outside of the domain $\Om$ of the equation. However
the similar techniques can be applied to our equations by the
stability property on viscosity subsolutions and supersolutions.

The key result of this section that is crucial for our regularity
theory is Theorem \ref{thm-5.4}, because we can apply it to
incremental quotients of viscosity solutions to fully nonlinear
integro-differential equations to get its $\rC^{1,\ap}$-estimates in
Section \ref{sec-10}.

In order to obtain a comparison principle for a nonlinear operator
$\cJ$, we need to impose a minimal ellipticity condition to our
collection $\fL$ of linear operators as follows (see also
\cite{CS}).

\begin{assumption}\label{ass-5.1} There is a constant $R_0\ge 1$ so that for
each $R>R_0$ and $\sm\in (1,2)$ there exists some $\dt=\dt(\sm,R)>0$
such that for any $\cL\in\fL$ we have that $\cL\vp>\dt\,\,\text{ on
$B_{R^{2-\sm}}$}$ where $\vp$ is a function given by
$\,\vp(x)=R^5\wedge |x|^2.$
\end{assumption}
Assumption \ref{ass-5.1} is enough for the comparison principle. We
note that Assumption \ref{ass-5.1} is very mild. In fact, we shall
show in the next lemma that the class $\fL_0$ satisfies Assumption
\ref{ass-5.1}. It just says that, given the particular function
$R^5\wedge|x|^2$, the value of the operator will be strictly
positive on $B_{R^{2-\sm}}$ but it does not require any uniform
estimate on how that happens.

Assumption \ref{ass-5.1} is pretty mild. Indeed, the following lemma
can be shown by simple computation. So we just state it without
detailed proof.

\begin{lemma}\label{lem-5.2} If $\,1<\sm<2$, then $\fL_0=\bigcup_{t\ge 1/2}\fL_t$
satisfies Assumption \ref{ass-5.1}.\end{lemma}

\begin{thm}\label{thm-5.3} Let $\cJ$ be elliptic with respect to $\fL$
in the sense of Definition \ref{def-3.1} where $\fL$ is some class satisfying
Assumption \ref{ass-5.1} and let $\Om\subset\BR^n$ be a bounded domain. If
$u\in\rB(\BR^n)$ is a viscosity subsolution to $\cJ u\ge f$ on
$\Om$, $v\in\rB(\BR^n)$ is a viscosity supersolution to $\cJ v\le f$
on $\Om$ and $u\le v$ on $\BR^n\s\Om$, then $u\le v$ on
$\Om$.\end{thm}

We obtain in Theorem \ref{thm-5.4} the result which allows functions
$u$ and $v$ to be discontinuous on $\BR^n\s\Om$ and is useful in
proving the comparison principle. Its proof can be done by a
nonsymmetric adaptation of the proof in \cite{CS}.

\begin{thm}\label{thm-5.4} Let $\cJ$ be elliptic with respect to some
class $\fL$ in the sense of Definition \ref{def-3.1} and let $\Om\subset\BR^n$
be a bounded domain. If $u\in\rB(\BR^n)$ is a viscosity subsolution
to $\cJ u\ge f$ on $\Om$ and $\,v\in\rB(\BR^n)$ is a viscosity
supersolution to $\cJ v\le g$ on $\Om$, then $\cM^+_{\fL}[u-v]\ge
f-g$ on $\Om$ in the viscosity sense.
\end{thm}

\begin{lemma}\label{lem-5.5} Let $\,\fL\,$ be a class of linear
integro-differential operators which satisfies Assumption \ref{ass-5.1}. If
$\,u\in\rB(\BR^n)$ is a viscosity subsolution to $\cM^+_{\fL} u\ge
0$ on $\Om$, then $\sup_{\Om}u\le\sup_{\BR^n\s\Om}u$.\end{lemma}

\pf Given $\sm\in(1,2)$, take a sufficiently large $R>0$ so that
$\Om\subset B_{R^{2-\sm}}$. For $\ep>0$ and $M\in\BR$, let
$\vp^{\ep}_M(x)=M+\f{\ep}{R^5}(R^5-R^5\wedge|x|^2).$ Then
$M\le\vp^{\ep}_M(x)\le M+\ep$ for any $x\in\BR^n$. Since $\cL\,1=0$
for all $\cL\in\fL$, by Assumption \ref{ass-5.1} there is a $\dt>0$
so that $\cM^+_{\fL}[\vp^{\ep}_M](x)\le-\ep\dt/R^5$ for any $x\in
B_{R^{2-\sm}}$. Then we can complete the proof by applying a similar
method as in \cite{CS}.\qed

{\it $[$Proof of Theorem \ref{thm-5.3}$]$} By Theorem \ref{thm-5.4},
we see that $\cM^+_{\fL}[u-v]\ge 0$ on $\Om$ in the viscosity sense.
Applying Lemma \ref{lem-5.5}, we obtain that
$\sup_{\Om}[u-v]\le\sup_{\BR^n\s\Om}[u-v].$ Hence this completes the
proof. \qed

\begin{remark} Once we obtain the comparison principle for viscosity
subsolutions and supersolutions which is semicontinuous on
$\overline\Om$, existence of the solutions of the Dirichlet problem
that we mentioned in the introduction follows from the Perron's
method \cite{I} as long as we can construct suitable barriers. For a
domain $\Om$ which has the exterior ball condition and prescribed
boundary data in $\BR^n\s\Om$ being continuous, the function $\Psi$
to be constructed in Section \ref{sec-7} can be used as the
barriers.\end{remark}
\section{ A nonlocal Alexandroff-Bakelman-Pucci estimate}\label{sec-6}

The Alexandroff-Bakelman-Pucci (A-B-P) estimate plays an important
role in  Krylov and Sofonov theory \cite{KS}, which is an essential
tool in the proof of Harnack inequality for linear uniformly
elliptic equations with measurable coefficients. In this section,
the influence of the gradient term has been addressed in our proof
of A-B-P estimate to which converges as $\sm$ is getting close to
$2$. In a later section, we shall use this nonlocal version of the
A-B-P estimte to prove H\"older estimates for $\sm$ close to $2$.

Let $u:\BR^n\to\BR$ be a function which is not positive outside the
ball $B_{1/2}$ and is upper semicontinuous on $\overline B_1$. We
consider its {\it concave envelope} $\Gm$ in $B_3$ defined as
$\Gm(x)=\inf\{p(x): p\in\Pi,\,p>u^+\,\,\text{in $B_2$}\}$ in $B_2$
and $0$ in $\BR^n\s B_2$, where $\Pi$ is the family of all the
hyperplanes in $\BR^n$. Also we denote the {\it contact set} of $u$
and $\Gm$ in $B_1$ by $\cC(u,\Gm,B_1)=\{y\in B_1:u(y)=\Gm(y)\}$.
\subsection{Key Lemma for nonlocal A-B-P estimate}
\begin{lemma}\label{lem-6.1} Let us assume $1<\sigma<2$. Let $u\le 0$ in
$\BR^n\s B_{1/2}$ and let $\Gm$ be its concave envelope in $B_2$. If
$\,u\in\rB(\BR^n)$ is a viscosity subsolution to $\cM^+_{\fL_0}
u=-f$ on $B_1$ where $f:\BR^n\to\BR$ is a function with $f>0$ on
$\cC(u,\Gm,B_1)$, then there exists some constant $C>0$ depending
only on $n,\ld$ and $\Ld$ $($but not on $\sm$$)$ such that for any
$x\in\cC(u,\Gm,B_1)$ and any $M>0$ there is some $k\in\BN\cup\{0\}$
such that
\begin{equation}\label{eq-6.1}
\bigl|\{y\in R_k(x):\underline{\mu}^-(u,x,y;\n\Gamma)\ge M
r_k^2\}\bigr|\le C\f{f(x)+|\n\Gm(x)|}{M}|R_k(x)|
\end{equation}
 where
$R_k(x)=B_{r_k}(x)\s B_{r_{k+1}}(x)$ for $r_k=\vr_0
2^{-\f{1}{2-\sm}-k}$ and $\vr_0=1/(16\sqrt n)$. Here, $\n\Gm(x)$
denotes any element of the superdifferential $\pa\Gm(x)$ of $\,\Gm$
at $x$. \end{lemma}

\rk We note that $\n\Gm(x)=\n u(x)$ for $x\in B_1$ if $\Gm$ and $u$
are differentiable at $x\in B_1$. In this case, $\pa\Gm(x)$ is a
singleton set with element $\n u(x)$.

{\it $[$Proof of Lemma \ref{lem-6.1}$]$} Take any
$x\in\cC(u,\Gm,B_1)$. Since $u$ can be touched by a hyperplane from
above at $x$, we see $\n\vp(x)=\n\Gm(x)$ for a
$\vp\in\rC^2_{\Om}(u;x)^+$. Thus it follows from Theorem
\ref{thm-3.4} that $\cM^+_{\fL_0} u(x;\n\Gm)$ is well-defined. We
observe that $\mu_t(u,x,y;\n\Gm) =u(x+y)-u(x)-(\n\Gm(x)\cdot
y)\,\chi_{B_t}(y)\le 0$ for any $y\in B_{1/2}$ and $t\ge 1/2$, by
the definition of concave envelope of $u$ in $B_2$. Since
$\mu_t^+(u,x,y;\n\Gm)\le|\n \Gm(x)|\,|y|$ for any $t\ge 1/2$, we
obtain that
\begin{equation*}\begin{split}\sup_{t\ge 1/2}\int_{\BR^n}\f{\Ld\mu_t^+(u,x,y;\n\Gm)}{|y|^{n+\sm}}\,dy
&\le\int_{\BR_n\s B_{1/2}}\f{\Ld |\n\Gm(x)|\,|y|}{|y|^{n+\sm}}\,dy
\fd c_0|\n\Gm(x)|.\end{split}\end{equation*} Then the constant
$c_0>0$ depending only on $n,\sm$ and $\Ld$ is finite for
$\sigma>1$. Thus by Theorem \ref{thm-3.4} we have that
\begin{equation*}\begin{split}&-f(x)\le\cM^+_{\fL_0} u(x;\n\Gm)\\
&=(2-\sm)\sup_{t\ge
1/2}\int_{\BR^n}\f{\Ld\mu_t^+(u,x,y;\n\Gm)-\ld\mu_t^-(u,x,y;\n\Gm)}{|y|^{n+\sm}}\,dy\\
&\le(2-\sm)\sup_{t\ge
1/2}\int_{\BR^n}\f{\Ld\mu_t^+(u,x,y;\n\Gm)}{|y|^{n+\sm}}\,dy
-(2-\sm)\inf_{t\ge
1/2}\int_{\BR^n}\f{\ld\mu_t^-(u,x,y;\n\Gm)}{|y|^{n+\sm}}\,dy\\
&\le(2-\sm)\int_{\BR_n\s B_{1/2}}\f{\Ld
|\n\Gm|\,|y|}{|y|^{n+\sm}}\,dy-(2-\sm)\int_{\BR^n}\f{\ld\underline{\mu}^-(u,x,y;\n\Gm)}{|y|^{n+\sm}}\,dy\\
&\le
(2-\sm)c_0|\n\Gm(x)|-(2-\sm)\int_{B_{r_0}(x)}\f{\ld\mu_1^-(u,x,y;\n\Gm)}{|y|^{n+\sm}}\,dy,
\end{split}\end{equation*} where $r_0=\vr_0 2^{-\f{1}{2-\sm}}$. Splitting the above integral
in the rings $R_k(x)$, we have that
\begin{equation*}
\begin{split} f(x)&\ge(2-\sm)\ld\sum_{k=0}^{\iy}\int_{R_k(x)}\f{\mu_1^-(u,x,y;\n\Gm)}{|y|^{n+\sm}}\,dy
-(2-\sm)c_0 |\n\Gm(x)|.
\end{split}\end{equation*}

Assume that the conclusion \eqref{eq-6.1} does not hold, i.e. for any $C>0$
there are some $x_0\in\cC(u,\Gm,B_1)$ and $M_0>0$ such that
\begin{equation}\label{eq-6.2}
\bigl|\{y\in R_k(x_0): \underline{\mu}^-(u,x_0,y;\n\Gm)\ge M_0
r_k^2\}\bigr|>C\,\f{f(x_0)+|\n\Gm
(x_0)|}{M_0}|R_k(x_0)|\end{equation} for all $k\in\BN\cup\{0\}$.
Since $-\mu_1\le\mu_1^-$ and $(2-\sm)\f{1}{1-2^{-(2-\sm)}}$ remains
bounded below for $\sm\in (1,2)$, it follows from \eqref{eq-6.2}
that
\begin{equation*}\begin{split} f(x_0)+ (2-\sigma) c_0
|\n\Gm(x_0)|&\ge(2-\sm)\ld\sum_{k=0}^{\iy}\int_{R_k(x_0)}\f{-\mu_1(u,x_0,y;\n\Gm)}{|y|^{n+\sm}}\,dy
\\
&\ge c(2-\sm)\sum_{k=0}^{\iy}M_0\f{r_k^2}{r_k^{\sm}} C\f{f(x_0)
+|\n\Gm(x_0)|}{M_0} \\
&\ge\f{c C(2-\sm)\rho_0^2}{1-2^{-(2-\sm)}}( f(x_0) + |\n\Gm(x_0)|)\\
&\ge c C (f(x_0) + |\n\Gm(x_0)|)\end{split}\end{equation*} for any $C>0$. Taking
$C$ large enough, we obtain a contradiction. Hence we are done. \qed

\noindent{\it Remark.} Lemma \ref{lem-6.1} would hold for any particular
choice of $\vr_0$ by modifying $C$ accordingly. The particular
choice $\vr_0=1/(16\sqrt n)$ is convenient for the proofs in Section
\ref{sec-7}.

\begin{lemma}\cite{CS}\label{lem-6.2} Let $\Gm$ be a concave function on $B_r(x)$
where $x\in\BR^n$ and let $h>0$. If $|\{y\in
S_r(x):\Gm(y)<\Gm(x)+(y-x)\cdot\n\Gm(x)-h\}|\le\ep\,|S_r(x)|$ for
any small $\ep>0$ where $S_r(x)=B_r(x)\s B_{r/2}(x)$, then we have
$\Gm(y)\ge\Gm(x)+(y-x)\cdot\n\Gm(x)-h$ for any $y\in
B_{r/2}(x)$.\end{lemma}

\begin{cor}\label{cor-6.3}
For any $\ep>0$, there is a constant $C>0$ such that for any
function $u$ with the same hypothesis as Lemma \ref{lem-6.1}, there
is some $r\in (0,\vr_0 2^{-\f{1}{2-\sm}})$ such that
\begin{equation*}\begin{split} &\f{|\{y\in S_r(x):u(y)<u(x)+(y-x)\cdot\n\Gm(x)-C( f(x) +|\n\Gm(x)|)
r^2\}|}{|S_r(x)|}\le\ep,\\&\int_{\overline Q}g_{\eta}(\n\Gm (y))\det
[D^2 \Gm(y)]^-\,dy\le C\sup_{y\in\overline Q}(1+\e^{-n} |f(y)|^n)\,
|Q|
\end{split}\end{equation*}
for any $\e>0$ and any cube $Q\subset B_{r/4}(x)$
with diameter $d$ such that $x\in\overline{Q}$ and $r/4<d<r/2$, where
$\vr_0=1/(16\sqrt n)$ and $g_{\eta}(z)=(|z|^{n/(n-1)}+\eta
^{n/(n-1)})^{1-n}$.\end{cor}

\pf The first part can be obtained by choosing $M=C (
f(x)+|\n\Gm(x)|)/\ep$ in Lemma \ref{lem-6.1}. Also the second part follows as
a consequence of Lemma \ref{lem-6.2} and concavity;
\begin{equation*}\begin{split}
\det[D^2 \Gm(x)]^-\le C ( f(x)+|\n\Gm(x)|)^n\le 4^n
C\frac{(1+\e^{-n}|f(x)|^n)}{g_{\eta}(\n\Gm(x))}.
\end{split}\end{equation*} Thus we have that
$g_{\eta}(\n\Gm(x))\det[D^2 \Gm(x)]^-\le 4^n C (1+\e^{-n}|f(x)|^n).$

Take any $y\in\cC(u,\Gm,B_1)\cap Q$ where $Q\subset B_{r/4}(x)$ is a
cube with diameter $d$ such that $x\in\overline Q$ and $r/4<d<r/2$.
Similarly to the above, we can obtain that
$g_{\eta}(\n\Gm(\cdot))\det[D^2 \Gm(\cdot)]^-\le 4^n C
(1+\e^{-n}|f(\cdot)|^n)$ $\aee$ on $Q$ because $\det[D^2
\Gm(\cdot)]^-=0$ $\aee$ on $Q\s\cC(u,\Gm,B_1)$ as in \cite{CC}.
Hence this implies the second part. \qed

\subsection{A nonlocal  A-B-P estimate}

We obtain a nonlocal version of Alexandroff-Bakelman-Pucci
estimate in the following theorem.

\begin{thm}\label{thm-6.4} Let $\,u$ and $\Gm$ be functions as in
Lemma \ref{lem-6.1}. Then there exist a finite family $\{Q_j\}_{j=1}^m$ of
open cubes $Q_j$ with diameters $d_j$ such that

$(a)$ any two cubes $Q_i$ and $Q_j$ do not intersect, $(b)$
$\cC(u,\Gm,B_1)\subset\bigcup_{j=1}^m\overline Q_j$,

$(c)$ $\cC(u,\Gm,B_1)\cap\overline Q_j\neq\phi$ for any $Q_j$, $(d)$
$d_j\le\vr_0 2^{-\f{1}{2-\sm}}$ where $\vr_0=1/(16\sqrt n)$,

$(e)$ $\int_{\overline Q_j}g_{\e}(\n\Gm (y))\det[D^2
\Gm(y)]^-\,dy\le C\sup_{\overline Q_j}(1+\e^{-n}|f|^n) |Q_j|$,

$(f)$ $|\{y\in 4\sqrt n\,Q_j:u(y)\ge\Gm(y)-C(\sup_{\overline Q_j} (f +|\n\Gm|))
d_j^2\}|\ge\gm |Q_j|$,

\noindent where the constants $C>0$ and $\gm>0$ depend on $n,\Ld$
and $\ld$ $($ but not on $\sm$$)$.\end{thm}

\pf In order to obtain such a family, we start by covering $B_1$
with a tiling of cubes of diameter $\vr_0 2^{-\f{1}{2-\sm}}$. Then
discard all those that do not intersect $\cC(u,\Gm,B_1)$. Whenever
a cube does not satisfy (e) and (f), we split it into $2^n$ cubes
of half diameter and discard those whose closure does not
intersect $\cC(u,\Gm,B_1)$. Now our goal is to prove that
eventually all cubes satisfy (e) and (f) and this process ends
after a finite number of steps.

Assume that the process does not finish in a finite number of steps.
Then we can have an infinite nested sequence of cubes. The
intersection of their closures will be a point $\hat x$. So we may
choose a sequence $\{x_k\}\subset\cC(u,\Gm,B_1)$ with
$\lim_{k\to\iy}x_k=\hat x$. Since $u(x_k)=\Gm(x_k)$ for all
$k\in\BN$, by the upper semicontinuity of $u$ on $\overline B_1$ we
have that $\Gm(\hat x)=\limsup_{k\to\iy}u(x_k)\le u(\hat x)$. Also
we have that $u(\hat x)\le\Gm(\hat x)$ because $u\le\Gm$ on $B_2$ by
the definition of the concave envelope $\Gm$ in $B_3$. Thus we
obtain that $u(\hat x)=\Gm(\hat x)$. We will now get a contradiction
by showing that eventually one of these cubes containing $\hat x$
will not split.

Take any $\ep>0$. By Corollary \ref{cor-6.3} there is some $r\in
(0,\vr_0 2^{-\f{1}{2-\sm}})$ such that
\begin{equation*}
\begin{split} &\f{|\{y\in
S_r(\hat x):u(y)<u(\hat x)+(y-\hat x)\cdot\n\Gm(\hat x)-C(f(\hat
x)+|\n\Gm(\hat x)|) r^2\}|}{|S_r(\hat x)|}\le\ep,\\&\int_{\overline
Q_j}g_{\eta}(\n\Gm (y))\det [D^2 \Gm(y)]^-\,dy\le
C\sup_{y\in\overline Q_j}(1+\e^{-n} |f(y)|^n)\, |Q_j|
\end{split}
\end{equation*}
for any $\e>0$ and a cube $Q_j\subset B_{r/4}(x)$ with diameter
$d_j$ such that $x\in\overline Q_j$ and $r/4<d_j<r/2$. So we easily
see that $\overline Q_j\subset B_{r/2}(\hat x)$ and $B_r(\hat
x)\subset 4\sqrt n\,Q_j$. We recall that $\Gm(y)\le u(\hat
x)+(y-\hat x)\cdot\n\Gm(\hat x)$ for any $y\in B_2$ because $\Gm$ is
concave on $B_2$ and $\Gm(\hat x)=u(\hat x)$. Since $d_j$ is
comparable to $r$, it thus follows that
\begin{equation*}\begin{split} &\bigl|\{y\in 4\sqrt n\,Q_j:u(y)\ge\Gm(y)-C(\sup_{\overline Q_j}(f+\n\Gm))d_j^2\}\bigr|\\
&\qquad\ge\bigl|\{y\in 4\sqrt n\,Q_j:u(y)\ge u(\hat x)+(y-\hat
x)\cdot\n\Gm(\hat x)-C(f(\hat x)+\n\Gm(\hat
x))r^2\}\bigr|\\&\qquad\ge(1-\ep)\bigl|S_r(\hat
x)\bigr|\ge\gm|Q_j|.\end{split}\end{equation*} Thus we proved (f).
Moreover, (e) holds for $Q_j$ because $\overline Q_j\subset
B_{r/2}(\hat x)$ and $B_r(\hat x)\subset 4\sqrt n\,Q_j$. Hence the
cube $Q_j$ would not split and the process must stop there. \qed

\begin{remark} Note that the upper bound for the diameters $\vr_0
2^{-\f{1}{2-\sm}}$ becomes very small as $\sm$ is getting close to
$2$. Adding $\sum_{j=1}^m |\n\Gm(Q_j)|$ and taking $\sm\to 2$, we
obtain the classical Alexandroff estimate as the limit of the
Riemann sums. For each $\sm>0$, we have that
$$\int_{\cC(u,\Gm,B_1) }g_{\e}(\n\Gm (y))\det[D^2 \Gm(y)]^-\,dy\le C\sum_j \sup_{y\in\overline Q_j}(1+\e^{-n}
|f(y)|^n)\, |Q_j|.$$ As $\sm\to 2$, the cube covering of
$\cC(u,\Gm,B_1)$ is getting close to the contact set
$\cC(u,\Gm,B_1)$ and the above becomes the following estimate
$($refer to Ch. 9.1 of \cite{GT}$)$;
$$\int_{B_{M_0}}g_{\e}(z)\,dz\le
C\int_{\cC(u,\Gm,B_1)}(1+\eta^{-n}|f(y)|^n)\,dy$$ for any $\e>0$,
where $M_0=\sup_{B_1} u^+$. Since $g_{\e}(z)\ge
2^{2-n}(|z|^n+\eta^n)^{-1}$, we have
$$\ln\biggl(\f{M_0^n}{\e^n}+1\biggr)\le
C\,\bigl(|B_1|+\e^{-n}\|f\|^n_{L^n(\cC(u,\Gm,B_1))}\bigr).$$ If we
set $\e=\|f\|_{L^n(\cC(u,\Gm,B_1))}$, then we obtain that
$\,\sup_{B_1} u^+\le C\,\|f\|_{L^n(\cC(u,\Gm,B_1))}.$ \end{remark}
\section{ Decay estimates of upper level sets}\label{sec-7}

In this section, we are going to apply the A-B-P estimate to get the
geometric decay rate of upper level sets for a nonnegative
subsolution $u$. To do this, we need a special function $\Psi$ so
that $\Psi-u$ meets the conditions of the A-B-P estimate. It is
based on the method used in \cite{CS}, but nontrivial computations
have been done to create positive terms so that it absorbs the
influence of the gradient term.
\subsection{Special Functions}

\begin{lemma}\label{lem-7.1} There exist some
$\sm^*\in (1,2)$ and $p>0$ such that the function $f(x)=2^p\wedge
|x|^{-p}$ is a subsolution to $\cM^-_{\fL_0} f(x)\ge 0$ for any
$\sm\in (\sm^*,2)$ and $x\in B_1^c$.\end{lemma}

\pf It is enough to show that there is some $\sm^*\in(1,2)$ so that
\begin{equation}\label{eq-7.1}
\cM^-_{\fL_0} f(x)\ge 0
\end{equation}
for $x=e_n=(0,\cdots,0,1)\in\BR^n$; for every other $x$ with
$|x|=1$, the above inequality follows by rotation. If $|x|\ge 1$,
then we consider $g(y)=|x|^p f(|x|y)$. Note that $g(y)=2^p|x|^p$ for
$|y|<1/(2|x|)$ and $|y|^{-p}$ for $|y|\ge 1/(2|x|)$. Then we can
derive that $g(y)\ge f(y)$ for any $y\in\BR^n$, $f(x/|x|)=g(x/|x|)$,
$\n f(x/|x|)=\n g(x/|x|)$ and
\begin{equation}\label{eq-7.2}
\n f\biggl(\f{x}{|x|}\biggr)\cdot y=|x|^p\bigl(\n
f(x)\cdot(|x|y)\bigr)
\end{equation}
for $y\in B_{1/2}$. Thus we see that
$\mu_t(g,x/|x|,y)\ge\mu_t(f,x/|x|,y)$ for any $t\ge 1/2$. We denote
by $\cK_0$ the class of all positive kernels satisfying
\eqref{eq-2.3} and \eqref{eq-2.4}. For $t\ge 1/2$ and $K\in\cK_0$,
we define the map $K_t$ given by $K_t(y)=t^{-n-\sm}K(y/t)$. Then it
is easy to check that the mapping $\cK_0\to\cK_0$ given by $K\mapsto
K_{|x|}$ is an isometry because $K\mapsto K_{1/|x|}$ is its inverse
mapping. Since $|\n f(x)|=p/|x|^{p+1}$, it follows from
\eqref{eq-7.2} and the change of variables that
\begin{equation*}\begin{split}\cM^-_{\fL_0} g(x/|x|)
&\le\cL^t g(x/|x|)\\
&=|x|^{\sm+p}\int_{\BR^n}\bigl[f(x+y)-f(x)-(\n f(x)\cdot y\bigr)\chi_{B_{t|x|}}(y)\bigr]K_{|x|}(y)\,dy\\
&\fd |x|^{\sm+p}\cL^{t|x|}_{|x|}f(x)
\end{split}\end{equation*} for any $\cL^t\in\fL_0$ and $\sm\in (\sm^*,2)$.
Since $t|x|\ge 1/2$ and the mapping $\cK_0\to\cK_0$ given by
$K\mapsto K_{|x|}$ is an isometry, by taking the infimum of both
sides on $\fL_0$ in the above inequality we obtain that
\begin{equation}\label{eq-7.3}
\cM^-_{\fL_0} g(x/|x|)\le|x|^{\sm+p}\cM^-_{\fL_0}
f(x)\end{equation}
 for any $\sm\in (\sm^*,2)$. By \eqref{eq-7.1}, \eqref{eq-7.2} and
 \eqref{eq-7.3}, we conclude that
\begin{equation*}\begin{split}\cM^-_{\fL_0} f(x)\ge\f{1}{|x|^{\sm+p}}\cM^-_{\fL_0}g(x/|x|)
\ge\f{1}{|x|^{\sm+p}}\cM^-_{\fL_0} f(x/|x|)\ge
0\end{split}\end{equation*} for any $\sm\in (\sm^*,2)$.

In order to prove \eqref{eq-7.1}, we use the following elementary
inequality that holds for any $a>b>0$ and $p>0$; $$(a+b)^{-p}\ge
a^{-p}\biggl(1-p\,\f{b}{a}+\f{p(p+1)}{2!}\bigl(\f{b}{a}\bigr)^2-\f{p(p+1)(p+2)}{3!}\bigl(\f{b}{a}\bigr)^3\biggr).$$
Using this inequality, we have that
\begin{equation}\label{eq-7.4}
\begin{split}\mu_t(f,e_n,y)&=|e_n+y|^{-p}-1+p\,y_n=(1+|y|^2+2 y_n)^{-p/2}-1+p\,y_n\\
&\ge-\bigl(\f{p}{2}+1\bigr)|y|^2+\f{p(p+2)}{2}
y_n|y|^2+\f{|y|^2}{(1+|y|^2)^{p/2+1}}\\
&+\f{p(p+2)}{2}\f{y_n^2}{(1+|y|^2)^{p/2+2}}-\f{p(p+2)(p+4)}{6}\f{y_n^3}{(1+|y|^2)^{p/2+3}}
\end{split}\end{equation}
for any $t\ge 1/2$ and $y\in B_{1/2}$. We now choose some
sufficiently large $p>0$ so that
\begin{equation}\label{eq-7.5}
\f{p(p+2)}{2(1+r^2)^{p/2+2}}\int_{S^{n-1}}\theta_n^2\,d\sigma(\theta)+\frac{\omega_n}{(1+r^2)^{p/2+1}}
-\bigl(\frac{p}{2}+1\bigr)\omega_n=\dt_0(r)>0\end{equation} for any
sufficiently small $r>0$, where $\om_n$ is the surface measure of
$S^{n-1}$. Since
$\int_{S^{n-1}}\theta_n\,d\sm(\theta)=\int_{S^{n-1}}\theta^3_n\,d\sm(\theta)=0$
and $\mu_t^-(f,e_n,y)\le 2^p+1+p$ for any $t\ge 1/2$ and $y\in
B_{1/2}$, it follows from \eqref{eq-3.1}, \eqref{eq-3.2},
\eqref{eq-7.4} and \eqref{eq-7.5} that
\begin{equation*}
\begin{split}&\cM^-_{\fL_0}
f(e_n)\\&\ge(2-\sm)\inf_{t\ge
1/2}\int_{\BR^n}\f{\ld\mu^+_t(f,e_n,y)}{|y|^{n+\sm}}\,dy
-(2-\sm)\sup_{t\ge
1/2}\int_{\BR^n}\f{\Ld\mu_t^-(f,e_n,y)}{|y|^{n+\sm}}\,dy\\
&\ge(2-\sm)\ld\,\inf_{t\ge
1/2}\int_{B_r}\f{\mu_t(f,e_n,y)}{|y|^{n+\sm}}\,dy
-(2-\sm)\Ld\,\sup_{t\ge
1/2}\int_{\BR^n\s B_r}\f{\mu_t^-(f,e_n,y)}{|y|^{n+\sm}}\,dy\\
&\ge(2-\sm)\biggl(\f{\ld\om_n\dt_0(r)}{2-\sm}
-(2^p+1+p)\Ld\int_{\BR^n\s B_r}\f{1}{|y|^{n+\sm}}\,dy\biggr)\\
&=\ld\om_n\dt_0(r)-(2^p+1+p)\Ld\,\om_n\f{2-\sm}{\sm}
\,r^{-\sm}\end{split}
\end{equation*}
for $r\in (0,1/2)$. Thus we may take some sufficiently small $r\in
(0,1/2)$ and take some $\sm^*\in (1,2)$ close enough to $2$ in the
above so that $\cM^-_{\fL_0} f(e_n)\ge 0$ for any $\sm\in
(\sm^*,2)$. Hence we complete the proof. \qed

\begin{cor}\label{cor-7.2} Given any $\sm_0\in(1,2)$, there exist some $\dt>0$ and $p>n$ such that the function
$f_{\dt}(x)=\dt^{-p}\wedge |x|^{-p}$ is a subsolution to
$\cM^-_{\fL_0} f_{\dt}(x)\ge 0$ for any $\sm\in (\sm_0,2)$, $x\in
B_1^c$.\end{cor}

\pf Let $\sm^*\in (1,2)$ be the number of Lemma \ref{lem-7.1}.
Without loss of generality, we may assume that $\sm_0<\sm^*$. Lemma
\ref{lem-7.1} implies that the result of this corollary always holds
for $\sm\in (\sm^*,2)$, when $\dt=1/2$. If $\dt<1/2$, then the
result still holds for $\sm\in (\sm^*,2)$ because
$\mu_t(f_{\dt},x,y)\ge\mu_t(f_{1/2},x,y)$ for any $y\in\BR^n$, $x\in
B_1^c$ and $t\ge 1/2$. We shall select $\dt\in (0,1-1/n)$ so small
that the result holds also for $\sm\in (\sm_0,\sm^*]$.

Now we let $x=e_n$ as in the proof of Lemma \ref{lem-7.1}. Assume
$\sm_0<\sm\le\sm^*$. Then we write
\begin{equation*}\begin{split}\cM^-_{\fL_0}
f_{\dt}(e_n)&=(2-\sm)\ld\int_{\BR^n}\f{{\underline\mu}^+(f_{\dt},e_n,y)}{|y|^{n+\sm}}\,dy
-(2-\sm)\Ld\int_{\BR^n}\f{{\overline\mu}^-(f_{\dt},e_n,y)}{|y|^{n+\sm}}\,dy\\
&\fd\cJ_1 f_{\dt}(e_n)+\cJ_2 f_{\dt}(e_n).\end{split}\end{equation*}
If we take some $\dt\in (0,1-1/n)$ small enough so that
$\mu_t^-(f_{\dt},e_n,y)=0$ for any $y\in B_{3/2}$ and $t\ge 1/2$,
from simple geometric observation it is easy to check that
$\mu_t^-(f_{\dt},e_n,y)\le 2^p +1+p|y|$ for any $y\in B_1^c$ and
$t\ge 1/2$. So we see that
$$-\cJ_2 f_{\dt}(e_n)=(2-\sm)\Ld\int_{|y|\ge
1}\f{{\overline\mu}^-(f_{\dt},e_n,y)}{|y|^{n+\sm}}\,dy\le(2-\sm_0)\Ld\int_{|y|\ge
1}\f{2^p +1+p|y|}{|y|^{n+\sm_0}}\,dy.$$ Since $\sm_0>1$, we have
that $\cJ_2 f_{\dt}(e_n)\ge-c_0$ for a constant $c_0>0$ depending on
$\sm_0,\Ld$ and the dimension $n$. On the other hand, since
$|y|<1+\dt$ and $-1+p y_n\ge-1+n(1-\dt)$ for any $y$ with
$\dt/2<|y+e_n|<\dt$, we have that
\begin{equation*}\begin{split}\cJ_1
f_{\dt}(e_n)&\ge(2-\sm)\ld\int_{\dt/2<|y+e_n|<\dt}\f{|e_n+y|^{-p}-1+n(1-\dt)}{|y|^{n+\sm}}\,dy\\
&\ge\f{(2-\sm)\ld}{(1+\dt)^{n+\sm}}\biggl(\int_{\dt/2<|y+e_n|<\dt}|e_n+y|^{-p}\,dy+(n-1-\dt n)\dt^n|B_1|(1-2^{-n})\biggr)\\
&\ge\f{(2-\sm)\ld}{(1+\dt)^{n+\sm}}\f{\om_n}{p-n}\dt^{n-p}(2^{p-n}-1).\end{split}\end{equation*}
If we select some $\dt\in(0,1-1/n)$ sufficiently small so that
$\cJ_1 f_{\dt}(e_n)>c_0,$ then we can complete the proof. \qed

\begin{lemma}\label{lem-7.3} Given any $\sm_0\in (1,2)$, there exists a function $\Psi\in\rB(\BR^n)$ such that

$(a)$ $\Psi$ is continuous on $\BR^n$, $(b)$ $\Psi=0$ on $B^c_{\sqrt
n}\,$,

$(c)$ $\Psi>2$ on $Q_1$, $(d)$ $\cM^-_{\fL_0}\Psi$ is continuous on
$B_{\sqrt n}\,$,

$(e)$ $\cM^-_{\fL_0}\Psi>-\psi$ on $\BR^n$ where $\psi$ is a
positive bounded function on $\BR^n$ which is supported in
$\overline B_{1/4}\,$, for any $\sm\in (\sm_0,2)$.\end{lemma}

\pf Let $\dt$ be the number of Corollary \ref{cor-7.2}. We consider
the function $\Psi$ given by $\Psi=0$ in $\BR^n\s B_{\sqrt n}$,
$c\bigl(|x|^{-p}-(\sqrt n)^{-p}\bigr)$ in $B_{\sqrt n}\s B_{\dt}$,
and $c\,P$ in $B_{\dt}$.
where $P$ is a quadratic
paraboloid chosen so that $\Psi$ is $\rC^{1,1}$ across $\pa
B_{\dt}$. We now choose the constant $c$ so that $\Psi(x)>2$ for
$x\in Q_1$ (recall that $Q_1\subset Q_2\subset B_{\sqrt n}\subset
B_{2\sqrt n}$). Since $\Psi\in\rC^{1,1}(B_{\sqrt n})$,
$\cM^-_{\fL_0}\Psi$ is continuous on $B_{\sqrt n}$. Also by
Corollary \ref{cor-7.2} we see that $\cM^-_{\fL_0}\Psi\ge 0$ on
$B^c_{1/4}$. Hence this completes the proof. \qed
\subsection{Key Lemma}
The main tool that shall be useful in proving H\"older estimates
is a lemma that connects a pointwise estimate with an estimate in
measure. The corresponding lemma in our context is the following.

\begin{lemma}\label{lem-7.4} Let $\sm_0\in (1,2)$ be given. If $\,\sm\in
(\sm_0,2)$, then there exist some constants $\vep_0>0$, $\nu\in
(0,1)$ and $M>1$ $($depending only on $\sm_0,\ld,\Ld$ and the
dimension $n$$)$ for which if $\,u\in\rB(\BR^n)$ is a viscosity
supersolution to $\cM^-_{\fL_0} u\le\vep_0$ on $B_{2\sqrt n}$ such
that $u\ge 0$ on $\BR^n$ and $\inf_{Q_1}u\le 1$, then $\,|\{u\le
M\}\cap Q_1|\ge\nu$.
\end{lemma}

\noindent{\it Remark.} (a) We denote by $Q_r(x)$ an open cube
$\{y\in\BR^n:|y-x|_{\iy}\le r/2\}$ and $Q_r=Q_r(0)$. If we set
$Q=Q_r(x)$, then we denote by $s Q=Q_{s r}(x)$ for $s>0$.

(b) If we assume that $0<\sm\le\sm^*<2$, then there is a simpler
proof of Lemma \ref{lem-7.3} using the ideas from \cite{S}. The
result here is more complicated as in \cite{CS} because we want to
get an estimate that remains uniform as $\sm\to 2$.

{\it $[$Proof of Lemma \ref{lem-7.4}$]$} We consider the function
$v\fd\Psi-u$ where $\Psi$ is the special function constructed in
Lemma \ref{lem-7.3}. Then we easily see that $v$ is upper
semicontinuous on $\overline B_{2\sqrt n}$ and $v$ is not positive
on $\BR^n\s B_{\sqrt n}$. Moreover, $v$ is a viscosity subsolution
to $\cM^+_{\fL_0} v\ge\cM^-_{\fL_0}\Psi -\cM^-_{\fL_0}
u\ge-\psi-\vep_0$ on $B_{2\sqrt n}$. So we want to apply Theorem
\ref{thm-6.4} (rescaled) to $v$. Let $\Gm$ be the concave envelope
of $v$ in $B_{4\sqrt n}$.

Since $\inf_{Q_1}u\le 1$, $\inf_{Q_1}\Psi>2$ and $Q_1\subset
B_{2\sqrt n}$, we easily see that $M_0\fd\sup_{B_{2\sqrt
n}}v=v(x_0)>0$ for some $x_0\in B_{2\sqrt n}$. We consider the
function $g$ whose graph is the cone in $\BR^n\times\BR$ with vertex
$(x_0,M_0)$ and base $\pa B_{6\sqrt n}(x_0)\times\{0\}$. For any
$\xi\in\BR^n$ with $|\xi|<M_0/6\sqrt n$, the hyperplane
$$H=\{(x,x_{n+1})\in\BR^n\times\BR:x_{n+1}=L(x)\fd
M_0+\xi\cdot(x-x_0)\}$$ is a supporting hyperplane for $g$ at $x_0$
in $B_{6\sqrt n}(x_0)$. Then $H$ has a parallel hyperplane $H'$
which is a supporting hyperplane for $v$ in $B_{4\sqrt n}$ at some
point $x_1\in B_{2\sqrt n}$. By the definition of concave envelope,
we see that $H'$ is also the hyperplane tangent to the graph of
$\Gm$ at $x_1$, so that $\xi=\n\Gm(x_1)$. This implies that
$B_{M_0/6\sqrt n}(0)\subset\n\Gm(B_{2\sqrt n})$. Thus we have that
\begin{equation}\label{eq-7.6}
C(n)\log \left(\frac{M_0^n}{\eta^n}\right)\leq \int_{\cC(u,\Gm,B_1) }g_{\e}(\n\Gm (y))
\det[D^2 \Gm(y)]^-\,dy,
\end{equation}
where $g_{\e}$ is the function given in Corollary \ref{cor-6.3}. We
also observe as shown in \cite{CC} that
\begin{equation}\label{eq-7.7}
\bigl|\n\Gm\bigl(B_{2\sqrt n}\s\cC(v,\Gm,B_{2\sqrt
n})\bigr)\bigr|=0.
\end{equation}
Let $\{Q_j\}$ be the finite family of cubes given by Theorem
\ref{thm-6.4} (rescaled on $B_{2\sqrt n}$). Then it follows from
\eqref{eq-7.6}, \eqref{eq-7.7} and Theorem \ref{thm-6.4} that
\begin{equation}\label{eq-7.8}\begin{split}\ln\biggl(\frac{[\sup_{B_{2\sqrt n}}v]^n}{\eta^n}+1\biggr)&\le
C\,\int_{\cC(u,\Gm,B_1) }g_{\e}(\n\Gm (y))\det[D^2\Gm(y)]^-\,dy\\
&\le C\biggl(\sum_j\sup_{\overline
Q_j}\bigl(1+\eta^{-n}(\psi+\vep_0)^n\bigr) |Q_j|\biggr)\\
&\le C\biggl(\sum_j|Q_j|+\e^{-n}\sum_j\sup_{\overline
Q_j}(\psi+\vep_0)^n |Q_j|\biggr).\end{split}
\end{equation}
If we set $\e=\bigl(\sum_j\sup_{\overline Q_j}(\psi+\vep_0)^n
|Q_j|\bigr)^{1/n}$ in \eqref{eq-7.8}, then we have that
\begin{equation}\label{eq-7.9}\begin{split}\sup_{B_{2\sqrt n}}v&\le C\biggl(\sum_j\sup_{\overline Q_j}(\psi+\vep_0)^n
|Q_j|\biggr)^{1/n}\le C\vep_0+C\biggl(\sum_j\bigl(\sup_{\overline
Q_j}\psi\bigr)^n |Q_j|\biggr)^{1/n}.\end{split}
\end{equation}
Since $\inf_{Q_1}u\le 1$ and $\inf_{Q_1}\Psi>2$, we see that
$\sup_{B_{2\sqrt n}}v>1$. If we choose $\e$ and $\vep_0$ small
enough, the above inequality \eqref{eq-7.9} implies that $1/2\le
C\bigl(\sum_j(\sup_{Q_j}\psi)^n |Q_j|\bigr)^{1/n}.$ We recall from
the proof of Lemma \ref{lem-7.3} that $\psi$ is supported on
$\overline B_{1/4}$ and bounded on $\BR^n$. Thus the above
inequality becomes $1/2\le C(\,\,\sum_{Q_j\cap B_{1/4}\neq\phi}
|Q_j|\,)^{1/n},$ which provides a lower bound for the sum of the
volumes of the cubes $Q_j$ intersecting $B_{1/4}$ as follows;
\begin{equation}\label{eq-7.10}\sum_{Q_j\cap
B_{1/4}\neq\phi} |Q_j|\ge c.
\end{equation}
Since $\diam(Q_j)\le\rho_0 2^{-\f{1}{2-\sm}}\le\rho_0$ for any
$\sm\in(\sm_0,2)$, the cube $4\sqrt n Q_j$ is contained in $B_{1/2}$
for any $Q_j$ with $Q_j\cap B_{1/4}\neq\phi$. Set
$M_1=\sup_{B_{1/2}}(\Psi-\Gm)$. Then by Theorem \ref{thm-6.4} we
have that
\begin{equation}\label{eq-7.11}\bigl|\{y\in 4\sqrt n
Q_j:v(x)\ge\Gm(y)-C d_j^2\}\bigr|\ge\gm |Q_j|\end{equation} and $C
d_j^2\le C\rho_0^2$. Then the family $\fF=\{4\sqrt n Q_j:Q_j\cap
B_{1/4}\neq\phi\}$ is an open covering of the union
$R\fd\bigcup\{\overline Q_j:Q_j\cap B_{1/4}\neq\phi\}.$ Now we may
take a subcovering of $\fF$ with finite overlapping number
(depending only on the dimension $n$) which covers the set $R$. Thus
it follows from \eqref{eq-7.10} and \eqref{eq-7.11} that
$\bigl|\{x\in B_{1/2}:v(x)\ge\Gm(x)-C\rho_0^2\}\bigr|\ge\nu.$ So we
have that $\,\bigl|\{x\in B_{1/2}:u(x)\le
M_1+C\rho_0^2\}\bigr|\ge\nu.$ Taking $M=M_1+C\rho_0^2>1$, we
conclude that $|\{u\le M\}\cap Q_1|\ge\nu$ because $B_{1/2}\subset
Q_1$. Hence we complete the proof. \qed

Let $Q_1$ be the unit cube. Then we split it into $2^n$ cubes of
half side. We do the same splitting step with each one of these
$2^n$ cubes and we continue this process. The cubes obtained in this
way are called {\it dyadic cubes.} If $Q$ is a dyadic cube different
from $Q_1$, then we say that $\Qtil$ is the {\it predecessor} of $Q$
if $Q$ is one of $2^n$ cubes obtained from splitting $\Qtil$.

\begin{lemma}\cite{CC}\label{lem-7.5} Let $A$ and $B$ be measurable subsets of $\BR^n$ with $A\subset
B\subset Q_1$. If $\,\dt\in (0,1)$ is some number such that $(a)$
$|A|\le\dt$ and $(b)$ $\Qtil\subset B$ for any dyadic cube $Q$ with
$|A\cap Q|>\dt|Q|$, then $|A|\le\dt |B|$.\end{lemma}

\subsection{Geometric Decay of Upper Level sets}
The following lemma is a consequence of Lemma \ref{lem-5.2}, Lemma \ref{lem-7.4} and
Lemma \ref{lem-7.5}.

\begin{lemma}\label{lem-7.6} Let $\vep_0>0$ be the constant in Lemma \ref{lem-7.4}.
If $u\in\rB(\BR^n)$ is a viscosity supersolution to
$\cM^-_{\fL_0}u\le\vep_0$ on $B_{2\sqrt n}$ such that $u\ge 0$ on
$\BR^n$ and $\inf_{Q_1}u\le 1$, then there exist universal constants
$C>0$ and $\vep_*>0$ such that $\bigl|\{u>t\}\cap Q_1\bigr|\le
C\,t^{-\vep_*}$ for any $t>0.$
\end{lemma}

\begin{remark}
 We note that $B_{1/2}\subset Q_1\subset Q_3\subset
B_{3\sqrt n/2}\subset B_{2\sqrt n}.$
\end{remark}

\pf Note that $B_{1/2}\subset Q_1\subset Q_3\subset B_{3\sqrt
n/2}\subset B_{2\sqrt n}.$ First, we shall prove that
\begin{equation}\label{eq-7.12}\bigl|\{u>M^k\}\cap
Q_1\bigr|\le(1-\nu)^k,\,\forall\,k\in\BN,\end{equation} where
$\nu>0$ and $M>1$ are the constants chosen as in Lemma
\ref{lem-7.4}.

If $k=1$, then it has been done in Lemma \ref{lem-7.4}. Assume the
result \eqref{eq-7.12} holds for $k-1$ ( $k\ge 2$ ), and let
$A=\{u>M^k\}\cap Q_1$ and $B=\{u>M^{k-1}\}\cap Q_1.$ If we can show
that $|A|\le(1-\nu)|B|$, then \eqref{eq-7.12} can be obtained for
$k$. To show this, we apply Lemma \ref{lem-7.6}. By Lemma
\ref{lem-7.4}, it is clear that $A\subset B\subset Q_1$ and
$|A|\le|\{u>M\}\cap Q_1|\le 1-\nu$. So it remain only to prove (b)
of Lemma \ref{lem-7.5}; that is, we need to show that if
$Q=Q_{2^{-i}}(x_0)$ is a dyadic cube satisfying
\begin{equation}\label{eq-7.13}|A\cap
Q|>(1-\nu)|Q|
\end{equation}
 then $\Qtil\subset B$. Indeed, we suppose
that $\Qtil\not\subset B$ and take $x_*\in\Qtil$ such that
\begin{equation}\label{eq-7.14}u(x_*)\le M^{k-1}.
\end{equation}
 We now consider the transformation
$x=x_0+2^{-i}y,\,\,y\in Q_1,\,\,x\in Q=Q_{2^{-i}}(x_0)$ and the
function $v(y)=u(x)/M^{k-1}$. If we can show that $v$ satisfies the
hypothesis of Lemma \ref{lem-7.4}, then we have that $\nu<|\{v(y)\le
M\}\cap Q_1|=2^{in}|\{u(x)\le M^k\}\cap Q|$, and thus $|Q\s A|>\nu
|Q|$ which contradicts \eqref{eq-7.13}.

To complete the proof, we consider once again the transformation
$$x=x_0+2^{-i}z,\,\,z\in B_{2\sqrt n},\,\,x\in B_{2\sqrt n/2^i}(x_0)\subset B_{2\sqrt n}$$ and the
function $v(z)=u(x)/M^{k-1}$. It now remains to show that $v$
satisfies the hypothesis of Lemma \ref{lem-7.4}. We now take any
$\vp\in\rC^2_{2\sqrt n}(v;z)^-$. If we set
$\psi=M^{k-1}\vp(2^i(\,\cdot\,-x_0))$, then we observe that
$$\vp\in\rC^2_{B_{2\sqrt n}}(v;z)^-\,\,\,\Leftrightarrow\,\,\,\psi\in
\rC^2_{B_{2\sqrt n/2^i(x_0)}}(u;x_0+2^{-i}z)^-.$$ If $K\in\cK_0$,
then we note that $K_i\in\cK_0$ where $K_i(y)=2^{i(n+\sm)}K(2^i y)$,
and moreover the mapping $\cK_0\to\cK_0$ given by $K\mapsto K_i$ is
an isometry. Since $B_{2\sqrt n/2^i}(x_0)\subset B_{2\sqrt n}$, we
have that
\begin{equation*}\begin{split}\cM^-_{\fL_0}v(z;\n\vp)&\le\cL^t
v(z;\n\vp)=\f{1}{2^{i\sm_0}M^{k-1}}\int_{\BR^n}\mu_{t
2^{-i}}(u,x_0+2^{-i}z,y;\n\psi)K_i(y)\,dy\\&\fd\f{1}{2^{i\sm_0}M^{k-1}}\cL_i^{t
2^{-i}} u(x_0+2^{-i}z;\n\psi)\end{split}\end{equation*} for any
$\cL^t\in\fL_t$ and any $t\ge 2^{i-1}$. Taking the infimum of the
right-hand side in the above inequality, we obtain that
$$\cM^-_{\fL_0}v(z;\n\vp)\le\f{1}{2^{i\sm_0}M^{k-1}}\cM^-_{\fL_0} u(x_0+2^{-i}z;\n\psi).$$
Thus $\cM^-_{\fL_0}v(z;\n\vp)\le\vep_0$ because
$\cM^-_{\fL_0}u\le\vep_0$ on $B_{2\sqrt n}\,$. By Theorem
\ref{thm-3.4}, we see that $\cM^-_{\fL_0}v\le\vep_0$ on $B_{2\sqrt
n}\,$. Also it is obvious that $v\ge 0$ on $\BR^n$ and we see from
\eqref{eq-7.14} that $\inf_{Q_1} v\le 1$. Finally the result follows
immediately from \eqref{eq-7.12} by taking $C=(1-\nu)^{-1}$ and
$\vep_*>0$ so that $1-\nu=M^{-\vep_*}$. Hence we complete the proof.
\qed

By a standard covering argument we obtain the following theorem.

\begin{thm}\label{thm-7.8} For any $\sm_0\in (1,2)$, let
$\sm\in(\sm_0,2)$ be given. If $u\in\rB(\BR^n)$ is a viscosity
supersolution to $\cM^-_{\fL_0} u\le\vep_0$ on $B_2$ such that $u\ge
0$ on $\BR^n$ and $u(0)\le 1$ where $\vep_0$ is the constant given
in Lemma \ref{lem-7.4}, then there are universal constants $C>0$ and
$\vep_*>0$ such that $\bigl|\{u>t\}\cap B_1\bigr|\le C\,t^{-\vep_*}$
for any $t>0.$\end{thm}

In contrast to symmetric cases, we can not obtain the following
theorem by rescaling the above theorem because our cases are not
scaling invariant. We note that Theorem \ref{thm-7.9} on $r\in
(0,1)$ shall be applied to obtain a Harnack inequality, and also
Theorem \ref{thm-7.9} on $r\in[1,2]$ will be used to prove an
interior $\rC^{1,\ap}$-regularity.

\begin{thm}\label{thm-7.9} For any $\sm_0\in (1,2)$, let
$\sm\in(\sm_0,2)$ be given, and let $x\in\BR^n$ and $r\in (0,2]$. If
$u\in\rB(\BR^n)$ is a viscosity supersolution to $\cM^-_{\fL_0} u\le
c_0$ on $B_{2r}(x)$ such that $u\ge 0$ on $\BR^n$, then there are
universal constants $\vep_*>0$ and $C>0$ such that
$$\bigl|\{u>t\}\cap B_r(x)\bigr|\le C\,r^n\bigl(u(x)+c_0\,r^{\sm}\bigr)^{\vep_*}t^{-\vep_*}\text{ for any $t>0$.}$$
\end{thm}

\pf Let $x\in\BR^n$ and set $v(z)=u(rz+x)/q$ for $z\in B_2$ where
$q=u(x)+c_0 r^{\sm}/\vep_0.$ Take any $\vp\in \rC^2_{B_2}(v;z)^-$.
If we set $\psi=q\,\vp((\,\cdot-x)/r)$, then we see that
$\psi\in\rC^2_{B_{2r}(x)}(u;rz+x)$. Thus by the change of variables
we have that
\begin{equation*}\begin{split}\cM^-_{\fL_0}v(z;\n\vp)&\le\cL^t
v(z;\n\vp)=\f{r^{\sm}}{q}\int_{\BR^n}\mu_{tr}(r z+x,y;\n\psi)\,K_r(y)\,dy\\
&\fd\f{r^{\sm}}{q}\cL^{tr}_r
u(rz+x;\n\psi)\end{split}\end{equation*} for any $\cL^t\in\fL_t$ and
any $t\ge 1/(2r)$, where $K_r(y)=r^{-n-\sm}K(y/r)$ for $r\in (0,2]$.
Taking the infimum of the right-hand side in the above inequality,
we get $$\cM^-_{\fL_0}v(z;\n\vp)\le\f{r^{\sm}}{q}\cM^-_{\fL_0}
u(rz+x;\n\psi)\le\vep_0.$$ Thus by Theorem 3.4 we have that
$\cM^-_{\fL_0}v\le\vep_0$ on $B_2$. Applying Theorem \ref{thm-7.8}
to the function $v$, we complete the proof. \qed
\section{Regularities}\label{sec-8.0}
\subsection{ Harnack inequality}\label{sec-8}

Harnack inequality plays an important role in analysis. In this
section, we obtain the Harnack inequality for integro-differential
equations whose associated kernel is not necessarily symmetric. Our
estimate depends only on a lower bound $\sm_0\in (1,2)$ for $\sm\in
(\sm_0,2)$ and also remains uniform as $\sm\to 2$. In this respect,
we can look upon this estimate as a generalization of Krylov-Safonov
Harnack inequality. The proof is an adaptation of the method used in
\cite{CS} to the case whose associated kernel is not necessarily
symmetric.

\begin{thm}\label{thm-8.1} For $\sm_0\in (1,2)$, let
$\sm_0<\sm<2$. If $\,u\in\rB(\BR^n)$ is a positive function such
that
$$\cM_{\fL_0}^- u\le C_0\,\,\,\text{ and }\,\,\,\cM_{\fL_0}^+ u\ge -C_0\,\,\text{
on $B_2$}$$ in the viscosity sense, then there is some constant
$C>0$ depending only on $\ld,\Ld,n$ and $\sm_0$ such that
$$\sup_{B_{1/2}} u\le C\,\bigl(\,\inf_{B_{1/2}}u+C_0\bigr).$$\end{thm}

\pf Let $\hat x\in B_{1/2}$ be a point so that
$\inf_{B_{1/2}}u=u(\hat x)$. Then it is enough to show that
$\sup_{B_{1/2}} u\le C\,\bigl(u(\hat x)+C_0\bigr).$ Without loss of
generality, we may assume that $u(\hat x)\le 1$ and $C_0=1$ by
dividing $u$ by $u(\hat x)+C_0$. Let $\vep_*>0$ be the number given
in Theorem \ref{thm-7.9} and let $\bt=n/\vep_*$. We now set
$s_0=\inf\{s>0:u(x)\le s(1-|x|)^{-\bt},\,\forall\,x\in B_1\}.$ Then
we see that $s_0>0$ because $u$ is positive on $\BR^n$. Also there
is some $\check x\in B_1$ such that $u(\check x)=s_0(1-|\check
x|)^{-\bt}=s_0\dd_0^{-\bt}$ where $\dd_0=\dd(\check x,\pa B_1)\le
1$.

To finish the proof, we have only to show that $s_0$ can not be too
large because $u(x)\le C_1(1-|x|)^{-\bt}\le C$ for any $x\in
B_{1/2}$ if $C_1>0$ is some constant with $s_0\le C_1$. Assume that
$s_0$ is very large. Then by Theorem \ref{thm-7.8} we have that
$$\bigl|\{u\ge u(\check x)/2\}\cap B_1\}\bigr|\le
\biggl|\f{2}{u(\check x)}\biggr|^{\vep_*}\le C
s_0^{-\vep_*}\dd_0^n.$$ Since $|B_{r}|=C\dd_0^n$ for $r=\dd_0/2<1$,
we easily obtain that
\begin{equation}\label{eq-8.1}\bigl|\{u\ge u(\check x)/2\}\cap
B_{r}(\check x)\}\bigr|\le\biggl|\f{2}{u(\check
x)}\biggr|^{\vep_*}\le C s_0^{-\vep_*}|B_{r}|.\end{equation} In
order to get a contradiction, we estimate $|\{u\le u(\check
x)/2\}\cap B_{\dt r}(\check x)|$ for some very small $\dt>0$ (to be
determined later). For any $x\in B_{2\dt r}(\check x)$, we have that
$u(x)\le s_0(\dd_0-\dt\dd_0)^{-\bt}\le u(\check x)(1-\dt)^{-\bt}$
for $\dt>0$ so that $(1-\dt)^{-\bt}$ is close to $1$.

We consider the function $v(x)=(1-\dt)^{-\bt}u(\check x)-u(x).$ Then
we see that $v\ge 0$ on $B_{2\dt r}(\check x)$, and also
$\cM_{\fL_0}^- v\le 1$ on $B_{\dt r}(\check x)$ because
$\cM_{\fL_0}^+ u\ge -1$ on $B_{\dt r}(\check x)$. We now want to
apply Theorem \ref{thm-7.9} to $v$. However $v$ is not positive on
$\BR^n$ but only on $B_{\dt r}(\check x)$. To apply Theorem
\ref{thm-7.9}, we consider $w=v^+$ instead of $v$. Since $w=v+v^-$,
we have that $\cM_{\fL_0}^- w\le\cM_{\fL_0}^- v+\cM_{\fL_0}^+ v^-\le
1+\cM_{\fL_0}^+ v^-$ on $B_{\dt r}(\check x)$. Since $v^-\equiv 0$
on $B_{2\dt r}(\check x)$, if $x\in B_{\dt r}(\check x)$ then we
have that $\mu_t(v^-,x,y;\n\vp)=v^-(x+y)$ for any $t\ge 1/2$, $y\in
B_{\dt r}(\check x)$ and $\vp\in\rC^2_{B_{\dt r}(\check
x)}(v^-;x)^+$. Take any $\vp\in\rC^2_{B_{\dt r}(\check x)}(v^-;x)^+$
and any $x\in B_{\dt r}(\check x)$. Since $x+B_{\dt r}\subset
B_{2\dt r}(\check x)$, we thus have that
\begin{equation}\label{eq-8.2}
\begin{split}&\cM_{\fL_0}^- w(x;\n\vp)\\&\qquad\le
1+(2-\sm)\int_{\BR^n}\f{\Ld\mu_t^+(v^-,x,y;\n\vp)-\ld\mu_t^-(v^-,x,y;\n\vp)}{|y|^{n+\sm}}\,dy\\
&\qquad\le 1+(2-\sm)\int_{\{y\in\BR^n:v(x+y)<0\}}\f{-\Ld \,v(x+y)}{|y|^{n+\sm}}\,dy\\
&\qquad\le 1+(2-\sm)\Ld\int_{\BR^n\s B_{\dt
r}}\f{\bigl(u(x+y)-(1-\dt)^{-\bt}u(\check
x)\bigr)_+}{|y|^{n+\sm}}\,dy.\end{split}\end{equation} Set
$h_c(x)=c(1-|x|^2)_+$ for $c>0$ and $c_1=\sup\{c>0:u(x)\ge
h_c(x),\forall x\in\BR^n\}.$ Then there is some $x_1\in B_1$ such
that $u(x_1)=c_1(1-|x_1|^2)$ and we see that $c_1\le 4/3$ because
$u(\hat x)\le 1$. Since $\n h_{c_1}(x)=-2 c_1 x$, we have that
\begin{equation}\label{eq-8.3}
\begin{split}(2-\sm)\int_{\BR^n}\f{\mu_t^-(u,x_1,y;\n
h_{c_1})}{|y|^{n+\sm}}\,dy&\le(2-\sm)\int_{\BR^n}\f{\mu_t^-(h_{c_1},x_1,y)}{|y|^{n+\sm}}\,dy\\
&\le\f{8(2-\sm_0)}{3}\int_{\BR^n}\f{1+|y|}{|y|^{n+\sm_0}}\,dy\le
C\end{split}\end{equation} for some constant $C>0$ which is
independent of $\sm$, and so we have that $$\Ld(2-\sm)\sup_{t\ge
1/2}\int_{\BR^n}\f{\mu_t^-(u,x_1,y;\n h_{c_1})}{|y|^{n+\sm}}\,dy\le
C.$$ Since $\cM_{\fL_0}^- u(x_1)\le 1$ on $B_2$, by \eqref{eq-8.3}
we have that
\begin{equation*}\begin{split} 1&\ge\cM_{\fL_0}^- u(x_1;\n h_{c_1})\ge\ld(2-\sm)\inf_{t\ge 1/2}\int_{\BR^n}\f{\mu_t^+(u,x_1,y;\n
h_{c_1})}{|y|^{n+\sm}}\,dy\\&\qquad\qquad\qquad\qquad-\Ld(2-\sm)\sup_{t\ge
1/2}\int_{\BR^n}\f{\mu_t^-(u,x_1,y;\n
h_{c_1})}{|y|^{n+\sm}}\,dy.\end{split}\end{equation*} Thus we obtain
that $\,\ds(2-\sm)\inf_{t\ge 1/2}\int_{\BR^n}\f{\mu_t^+(u,x_1,y;\n
h_{c_1})}{|y|^{n+\sm}}\,dy\le C\,$ for a constant $C>0$ which is
independent of $\sm$, and so there is some $t_1\ge 1/2$ such that
\begin{equation}\label{eq-8.4}(2-\sm)\int_{\BR^n}\f{\mu_{t_1}^+(u,x_1,y;\n
h_{c_1})}{|y|^{n+\sm}}\,dy\le C.
\end{equation}
Since $\mu_{t}(u,x_1,y;\n h_{c_1})\ge u(x_1 +y)-4/3-8t/3$ for any
$t\ge 1/2$ and $y\in B_t$, and $$(u(x_1
+y)-4/3-8t/3)_+=0\,\,\,\text{ for any }\,\,t\ge\f{3}{8}\sup_{y\in
\BR^n}\bigl[u(x_1+y)-4/3\bigr],$$ we may assume that $t_1\ge 1/2$
must be finite. Since $\mu^+_{t_1}(u,x_1,y;\n h_{c_1})\ge (u(x_1
+y)-4/3-8t_1/3)_+$ for any $y\in\BR^n$, by \eqref{eq-8.4} we obtain
that
\begin{equation}\begin{split}\label{eq-8.5}&(2-\sm)\int_{\BR^n}\f{\bigl(u(x_1+y)-4/3-8t_1/3\bigr)_+}{|y|^{n+\sm}}\,dy\\
&\qquad\qquad\le(2-\sm)\int_{\BR^n}\f{\mu_{t_1}^+(u,x_1,y;\n
h_{c_1})}{|y|^{n+\sm}}\,dy\le C\,. \end{split}\end{equation} We now
may assume that $(1-\dt)^{-\bt}u(\check
x)=(1-\dt)^{-\bt}s_0(1-|\check x|)^{-\bt}\ge 4/3+8 t_1/3$ because
$s_0$ was very large and $(1-\dt)^{-\bt}$ was close to $1$. Since
$\dt r<1$, by \eqref{eq-8.5} and the change of variables we have
that
\begin{equation*}\begin{split}&(2-\sm)\Ld\int_{B^c_{\dt
r}}\f{\bigl(u(x+y)-(1-\dt)^{-\bt}u(\check x)\bigr)_+}{|y|^{n+\sm}}\,dy\\
&\le(2-\sm)\Ld\int_{B^c_{\dt r}\cap
B_{100}}\f{\bigl(u(x_1+y+x-x_1)-(1-\dt)^{-\bt}u(\check
x)\bigr)_+}{|y+x-x_1|^{n+\sm}}\,
\f{|y+x-x_1|^{n+\sm}}{|y|^{n+\sm}}\,dy\\
&+(2-\sm)\Ld\int_{B^c_{\dt r}\cap
B_{100}^c}\f{\bigl(u(x_1+y+x-x_1)-(1-\dt)^{-\bt}u(\check
x)\bigr)_+}{|y+x-x_1|^{n+\sm}}\,
\f{|y+x-x_1|^{n+\sm}}{|y|^{n+\sm}}\,dy\\&\le C\bigl((\dt
r)^{-n-\sm}+1\bigr)(2-\sm)\Ld\int_{\BR^n}\f{\bigl(u(x_1+y)-4/3-8
t_1/3\bigr)_+}{|y|^{n+\sm}}\,dx\le C(\dt
r)^{-n-\sm}\end{split}\end{equation*} for any $x\in B_{\dt r}(\check
x)$. Thus by \eqref{eq-8.2} and Theorem \ref{thm-3.4} we obtain that
$$\cM_{\fL_0}^- w(x)\le C(\dt r)^{-n-\sm}\,\,\text{ on $B_{\dt
r}(\check x)$ }.$$ Since $u(\check x)=s_0\dd_0^{-\bt}=2^{-\bt}s_0
r^{-\bt}$ and $\bt\vep_*=n$, applying Theorem \ref{thm-7.9} we have
\begin{equation*}\begin{split}&\bigl|\{u\le u(\check x)/2\}\cap B_{\dt r/2}(\check x)\bigr|
=\bigl|\{w\ge u(\check x)((1-\dt)^{-\bt}-1/2)\}
\cap B_{\dt r/2}(\check x)\bigr|\\
&\qquad\le C(\dt r)^n\bigl[((1-\dt)^{-\bt}-1)u(\check x)+C(\dt
r)^{-n-\sm}(\dt
r)^{\sm}\bigr]^{\vep_*}\bigl[u(\check x)((1-\dt)^{-\bt}-1/2)\bigr]^{-\vep_*}\\
&\qquad\le C(\dt r)^n\bigl[((1-\dt)^{-\bt}-1)^{\vep_*}+\dt
^{-n\vep_*}s_0^{-\vep_*}\bigr].\end{split}\end{equation*} We now
choose $\dt>0$ so small enough that $C(\dt
r)^n((1-\dt)^{-\bt}-1)^{\vep_*}\le |B_{\dt r/2}(\check x)|/4.$ Since
$\dt$ was chosen independently of $s_0$, if $s_0$ is large enough
for such fixed $\dt$ then we get that $C(\dt r)^n\dt
^{-n\vep_*}s_0^{-\vep_*}\le |B_{\dt r/2}(\check x)|/4.$ Therefore we
obtain that $\bigl|\{u\le u(\check x)/2\}\cap B_{\dt r/2}(\check
x)\bigr|\le |B_{\dt r/2}(\check x)|/2.$ Thus we conclude that
\begin{equation*}\begin{split}\bigl|\{u\ge u(\check x)/2\}\cap B_r(\check x)\bigr|&\ge\bigl|\{u\ge u(\check x)/2\}\cap
B_{\dt r/2}(\check x)\bigr|\ge\bigl|\{u>u(\check x)/2\}\cap B_{\dt
r/2}(\check x)\bigr|\\&\ge\bigl|B_{\dt r/2}(\check
x)\bigr|-\bigl|B_{\dt r/2}(\check x)\bigr|/2=\bigl|B_{\dt
r/2}(\check x)\bigr|/2=C |B_r|,\end{split}\end{equation*} which
contradicts \eqref{eq-8.1} if $s_0$ is large enough. Thus we
complete the proof. \qed
\subsection{ H\"older estimates}\label{sec-9}

The purpose of this section is to prove the following H\"older
regularity result (see Theorem \ref{thm-9.2}).  Before doing this,
we obtain a technical lemma which shall be useful in proving Theorem
\ref{thm-9.2}.

\begin{lemma}\label{lem-9.1} For any $\sm_0\in (1,2)$, let $\sm\in
(\sm_0,2)$ be given. If $\,u$ is a bounded function with $|u|\le
1/2$ on $\BR^n$ such that
$$\cM^+_{\fL_0} u\ge-\vep_0\,\,\text{ and }\,\,\cM^-_{\fL_0} u\le\vep_0\,\,\text{
on $B_1$ }$$ in the viscosity sense where $\vep_0>0$ is some
sufficiently small constant, then there is some universal constant
$\ap>0$ $($depending only on $\ld,\Ld,n$ and $\sm_0$$)$ such that
$u\in\rC^{\ap}$ at the origin. More precisely, $|u(x)-u(0)|\le
C\,|x|^{\ap}$ for some universal constant $C>0$ depending only on
$\ap$.\end{lemma}

We shall prove the following theorem using only Theorem
\ref{thm-7.9}. Theorem \ref{thm-9.2} easily follows from Lemma
\ref{lem-9.1} by a simple rescaling argument. Thus we have only to
prove Lemma \ref{lem-9.1} to establish Theorem \ref{thm-9.2}.

\begin{thm}\label{thm-9.2} For any $\sm_0\in (1,2)$, let $\sm\in
(\sm_0,2)$ be given. If $u$ is a bounded function on $\BR^n$ such
that $\cM^+_{\fL_0} u\ge -c_0\,\text{ and }\,\cM^-_{\fL_0} u\le
c_0\,\text{ on $B_1$ }$ in the viscosity sense, then there is a
constant $\ap>0$ $($depending only on $\ld,\Ld,n$ and $\sm_0$$)$
such that $$\|u\|_{\rC^{\ap}(B_{1/2})}\le
C\bigl(\,\|u\|_{L^{\iy}(\BR^n)}+c_0\bigr)$$ where $C>0$ is some
universal constant depending only on $\ap$.\end{thm}

{\it $[$Proof of Lemma \ref{lem-9.1}$]$} Take any $\ap\in (0,\sm_0)$
and choose an $N$ so large that
\begin{equation}\label{eq-9.1}\f{2^{\sm_0+1}(2-\sm_0)\Ld\om_n}{\sm_0-\ap}\,2^{-(\sm_0-\ap)N}\le\vep_0/2\,\,\,\,\text{ and
}\,\,\,\, 2^{1-\sm_0 N} 2^{-k(\sm_0 -\ap)N}\le 1/2.\end{equation}
Then it suffices to show that there exist a nondecreasing sequence
$\{m_k\}_{k\in\BN\cup\{0\}}$ and a nonincreasing sequence
$\{M_k\}_{k\in\BN\cup\{0\}}$ such that $m_k\le u\le M_k$ in
$B_{2^{-kN}}$ and $M_k-m_k=2^{-\ap kN},$ so that the theorem holds
with $C=2^{\ap N}$; for, if $2^{-(k+1)N}\le|x|\le 2^{-kN}$ for
$k\in\BN\cup\{0\}$, then
$$\bigl|u(x)-u(0)\bigr|\le\f{M_k-m_k}{2^{-\ap(k+1)N}}\cdot
2^{-\ap(k+1)N}\le 2^{\ap N} |x|^{\ap}.$$ We construct $m_k$ and
$M_k$ by induction. For $k\le 0$, we can take $m_k=\inf_{\BR^n} u$
and $M_k=m_k+2^{-\ap kN}$ because $\osc_{\BR^n} u\le 1$. Assume that
we have the sequences up to $m_k$ and $M_k$ for $k\ge 1$. Then we
want to show that we can continue the sequences by finding $m_{k+1}$
and $M_{k+1}$. Fix any $z\in B_{1/\tau}$ where $\tau=2^{-(k+1)N}$.
Take any $\vp\in \rC^2_{B_{1/\tau}}(v;z)^-$ where
$v=\ds\f{u(\tau\cdot)-m_k}{(M_k-m_k)/2}.$ If we set
$\psi=m_k+\f{M_k-m_k}{2}\vp(\cdot/\tau)$, then we see that
$\psi\in\rC^2_{B_1}(u;\tau z)^-$ and moreover we have that
\begin{equation}\label{eq-9.2}\n\psi(\tau
z)=\f{M_k-m_k}{2\tau}\,\n\vp(z).\end{equation}

In $B_{2^{-(k+1)N}}$, either $u>(M_k+m_k)/2$ in at least half of the
points (in measure) or $u\le (M_k+m_k)/2$ in at least half of the
points. First, we deal with the case $$|\{u>(M_k+m_k)/2\}\cap
B_{2^{-(k+1)N}}|\ge |B_{2^{-(k+1)N}}|/2.$$ Then we note that $v\ge
0$ on $B_{2^N}$ and $|\{v>1\}\cap B_1|\ge |B_1|/2$. Then we have
\begin{equation*}\begin{split}\cM^-_{\fL_0} v(z;\n\vp)&\le\cL^t
v(z;\n\vp)=\f{2\,\tau^{\sm}}{M_k-m_k}\int_{\BR^n}\mu_{t\tau}(u,\tau
z,y;\n\psi)\,K_{\tau}(y)\,dy\\
&\fd\f{2\,\tau^{\sm}}{M_k-m_k}\cL^{t\tau}_{\tau} u(\tau
z;\n\psi)\end{split}\end{equation*} for any $\cL^t\in\fL_t$ and any
$t\ge 2^{(k+1)N-1}=1/(2\tau)$. Taking the infimum of the right-hand
side in the above inequality, we obtain that
$$\cM^-_{\fL_0}v(z;\n\vp)\le\f{2\,\tau^{\sm}}{M_k-m_k}\cM^-_{\fL_0} u(\tau z;\n\psi).$$ Since
$\cM^-_{\fL_0} u\le\vep_0$ on $B_1$ in the viscosity sense, we have
that
\begin{equation*}\cM^-_{\fL_0}v(z;\n\vp)
\le \f{2\,\tau^{\sm}}{M_k-m_k}\cM^-_{\fL_0} u(\tau z;\n\psi)\le
2^{1-\sm_0 N}2^{-k(\sm_0-\ap)N}\vep_0\le\vep_0/2.\end{equation*}
Thus this implies that
\begin{equation}\label{eq-9.3}\cM^-_{\fL_0}v\le\vep_0/2\,\,\text{ on $B_{2^{(k+1)N}}$.
}
\end{equation}
It also follows from the inductive hypothesis that if $2^{jN}\le
|x|\le 2^{(j+1)N}$ then
\begin{equation}\label{eq-9.4}
\begin{split}
v(x)&\ge\f{m_{k-j}-M_{k-j}+M_k-m_k}{(M_k-m_k)/2}\ge 2(1-|x|^{\ap}),\\
v(x)&\le\f{M_{k-j}-m_{k-j}+m_{k-j}-m_k}{(M_k-m_k)/2}\le
2|x|^{\ap}.\end{split}\end{equation} for any $j\in\BN\cup\{0\}$.
That is to say, $-2(|x|^{\ap}-1)\le v(x)\le 2|x|^{\ap}$ outside
$B_1$.

If we set $w(x)=\max\{v(x),0\}$, then we have that $\cM^-_{\fL_0}
w\le 4\vep_0$ on $B_{2^{N-1}}$; for, if $x\in B_{2^{N-1}}$ and $y\in
B_{2^{N-1}}$, then $\mu_t(v^-,x,y)=0$ for any $t\ge 1/2$ because
$v\ge 0$ on $B_{2^N}$; if $x\in B_{2^{N-1}}$ and $y\in
B^c_{2^{N-1}}$, then $\mu_t(v^-,x,y)=v^-(x+y)$ for any $t\ge 1/2$
because $v\ge 0$ on $B_{2^N}$. Since $w=v+v^-$, we see that
$\cM^-_{\fL_0} w\le\cM^-_{\fL_0} v+\cM^+_{\fL_0} v^-$. Thus it
follows from \eqref{eq-9.1} and \eqref{eq-9.4} that if $x\in
B_{2^{N-1}}$ is given, then
\begin{equation*}\begin{split}\cL^t v^-(x)&\le\int_{\BR^n}\mu_t(v^-,x,y)K(y)\,dy\le
(2-\sm)\Ld\int_{|y|\ge 2^{N-1}}\f{2\,|x+y|^{\ap}}{|y|^{n+\sm}}\,dy\\
&\le(2-\sm_0)\Ld\int_{|y|\ge
2^{N-1}}\f{2^{\ap+1}}{|y|^{n+\sm_0-\ap}}\,dy
=\f{2^{\sm_0+1}(2-\sm_0)\Ld\om_n}{\sm_0-\ap}\,2^{-(\sm_0-\ap)N}\le\vep_0/2\end{split}\end{equation*}
for any $\cL^t\in\fL_0$, whenever $\ap\in (0,\sm_0)$ and
$\sm\in(\sm_0,2)$. So we have that $\cM^+_{\fL_0} v^-\le\vep_0/2$ on
$B_{2^{N-1}}$. Thus by \eqref{eq-9.3} we conclude that
$\cM^-_{\fL_0} w\le \vep_0$ on $B_{2^{N-1}}$ where $\ap\in
(0,\sm_0)$.

Take any point $x\in B_1$. Since $B_1\subset B_2(x)\subset
B_4(x)\subset B_{2^{N-1}}$, we can apply Theorem \ref{thm-7.9} on
$B_2(x)$ to obtain that $2^n
C(w(x)+2^{\sm}\vep_0)^{\vep_*}\ge\bigl|\{w>1\}\cap
B_2(x)\bigr|\ge\bigl|\{v>1\}\cap B_1\bigr|\ge\f{1}{2}\,|B_1|.$ Thus
we have that $\ds
w(x)\ge\biggl(\f{|B_1|}{2^{n+1}C}\biggr)^{1/\vep_*}-4\,\vep_0$ for
any $x\in B_1$. If we choose $\vep_0$ sufficiently small, then this
implies that $w\ge\vartheta$ on $B_1$ for some $\vartheta>0$. If we
set $M_{k+1}=M_k$ and $m_{k+1}=m_k+\vartheta(M_k-m_k)/2$, then we
have that $m_{k+1}\le u\le M_{k+1}$ in $B_{2^{-(k+1)N}}$. Moreover,
$M_{k+1}-m_{k+1}=(1-\vartheta/2)2^{-\ap kN}$. Then we may choose
some small $\ap>0$ and $\vartheta>0$ so that $1-\vartheta/2=2^{-\ap
N}$, so that we obtain that $M_{k+1}-m_{k+1}=2^{-\ap(k+1)N}$.

On the other hand, if we deal with the second case
$$|\{u\le(M_k+m_k)/2\}\cap B_{2^{-(k+1)N}}|\ge
|B_{2^{-(k+1)N}}|/2,$$ then we consider the function $\ds
v(x)=\f{M_k-u(2^{-(k+1)N}x)}{(M_k-m_k)/2}$ and repeat in the same
way by using $\cM^+_{\fL_0} u\ge-\vep_0$.\qed

\subsection{ $\rC^{1,\ap}$-estimates}\label{sec-10}

In this section, we prove an interior $\rC^{1,\ap}$-regularity
result for viscosity solutions to a general class of fully nonlinear
integro-differential equations. The key idea of proof is to apply
the H\"older estimates of Theorem \ref{thm-9.2} to incremental
quotients of the solution. There being no uniform bound in $L^{\iy}$
for the incremental quotients outside the domain may cause a
technical difficulty since we are dealing with nonlocal equations.
We shall solve it by assuming some extra regularity of the family of
linear integro-differential operators $\cL^t$. The extra assumption,
added to the growth condition \eqref{eq-2.3} for the kernel, is a
modulus of continuity of $K$ in measure so that far away
oscillations tend to cancel out.

We consider the class $\fL^1_0$ consisting of the operators
$\cL^t\in\fL_0$ associated with kernels $K$ for which \eqref{eq-2.3}
holds and there exists some $\vr_1>0$ such that
\begin{equation}\label{eq-10.1}\sup_{h\in B_{\vr_1/2}}\int_{\BR^n\s B_{\vr_1}}\f{|K(y)-K(y-h)|}{|h|}\,dy\le
C.\end{equation}

If $K$ is a radial function satisfying \eqref{eq-2.3}, then it is
interesting that the condition \eqref{eq-10.1} is not required.
Indeed, if $K(y)=(2-\sm)A/|y|^{n+\sm}$ for $\ld\le A\le\Ld$, then it
follows from the mean value theorem and Schwartz inequality that
\begin{equation*}\begin{split}\sup_{h\in B_{\vr_1/2}}\int_{\BR^n\s
B_{\vr_1}}\f{|K(y)-K(y-h)|}{|h|}\,dy&=(2-\sm)\Ld(n+\sm)2^{n+\sm+1}\f{\om_n}{\sm+1}\vr_1^{-1-\sm}\le
C\,.\end{split}\end{equation*} for any $h\in B_{\vr_1/2}$ and
$y\in\BR^n\s B_{\vr_1}$, because $|y|\ge 2|h|$ for such $h,y$ and
$|y-\tau h|\ge |y|-|h|\ge |y|-|y|/2=|y|/2$ for $\tau\in [0,1]$.

In the following theorem, we shall furnish interior
$\rC^{1,\ap}$-estimates for fully nonlinear elliptic equations
associated with a class of linear integro-differential operators.

\begin{thm}\label{thm-10.1} For any $\sm_0\in (1,2)$, let $\sm\in
(\sm_0,2)$ be given. Then there is some $\vr_1>0$ $($depending on
$\ld,\Ld,\sm_0$ and the dimension $n$$)$ so that if $\,\cI^{\pm}$ is
a nonlocal elliptic operator with respect to $\fL^1_0$ in the sense
of Definition \ref{def-3.1} and $u\in\rB(\BR^n)$ is a viscosity
solution to $\cI^{\pm} u=0$ on $B_1$, then there is a universal
constant $\ap>0$ $($depending only on $\ld,\Ld,\sm_0$ and the
dimension $n$$)$ such that
$$\|u\|_{\rC^{1,\ap}(B_{1/2})}\le
C\bigl(\,\|u\|_{L^{\iy}(\BR^n)}+|\cI^{\pm} 0|\,\bigr)$$ for some
constant $C>0$ depending on $\ld,\Ld,\sm_0,n$ and the constant given
in $($$10.1$$)$ $($where we denote by $\cI^{\pm} 0$ the value we
obtain when we apply $\cI^{\pm}$ to the constant function that is
equal to zero$)$.\end{thm}

\pf Since $\cI^{\pm} u=0$ on $B_1$, it follows from Lemma
\ref{lem-3.2} that $\cM^+_{\fL_0} u\ge\cI^{\pm} u-\cI^{\pm}
0=-\cI^{\pm} 0\ge-|\cI^{\pm} 0|$ on $B_1$. Similarly we have that
$\cM^-_{\fL_0} u\le|\cI^{\pm} 0|$ on $B_1$. Thus by Theorem
\ref{thm-9.2} we see that $u\in\rC^{\ap}(B_{1-\dt})$ for any
$\dt\in(0,1)$ and $\|u\|_{\rC^{\ap}(B_{1-\dt})}\le
C(\,\|u\|_{L^{\iy}(\BR^n)}+|\cI^{\pm} 0|\,).$ Now we want to improve
the obtained regularity iteratively by applying Theorem
\ref{thm-9.2} again until we obtain Lipschitz regularity in a finite
number of steps.

Assume that we have proven that $u\in\rC^{\bt}(B_r)$ for some
$\bt\in(0,1]$ and $r\in (0,1)$. Then we apply Theorem \ref{thm-9.2}
for the difference quotient $w^h=(\tau_h u-u)/|h|^{\bt}$ where
$\tau_h$ is a translation operator given by $\tau_h u(x)=u(x+h)$ for
$h\in\BR^n$. Since $\cM^+_{\fL_0} u=-\cM^-_{\fL_0}[-u]$, we see from
Theorem \ref{thm-5.4} that $\cM^+_{\fL_0} w^h\ge 0$ and
$\cM^-_{\fL_0} w^h\le 0$ on $B_r$ for any $h\in (0,1-r)$. Since
$u\in\rC^{\bt}(B_r)$, we see that the family
$\{w^h\}_{|h|\in(0,1-r)}$ is uniformly bounded on $B_r$. But the
functions $w^h$ is not uniformly bounded outside the ball $B_r$, and
so we can not apply Theorem \ref{thm-9.2} directly. However we
observe from \eqref{eq-10.1} that $w^h$ has oscillations that give
cancelations in the integrals. Let $\phi$ be a smooth cutoff
function supported in $B_r$ such that $\phi\equiv 1$ in
$B_{r-\dt/4}$ where $\dt>0$ is some small positive number which
shall be determined later. We write $w^h=w_1^h+w_2^h$ where
$w_1^h=\phi\,w^h$ and $w_2^h=(1-\phi)\,w^h$. Take any $x\in
B_{r-\dt/2}$ and $|h|<\dt/16$. Then
$(1-\phi(x))u(x)=(1-\phi(x))\tau_h u(x)=0$ and $w^h(x)=w_1^h(x)$. We
prove that $w_1^h\in\rC^{\ap+\bt}(B_{r-\dt})$ for some $\ap>0$ with
$\ap+\bt>1$. We note that
\begin{equation}\label{eq-10.2}
\begin{split}\cM^+_{\fL_0} w_1^h\ge\cM^+_{\fL^1_0} w_1^h&=\cM^+_{\fL^1_0} [w^h-w_2^h]\ge
0-\cM^+_{\fL^1_0} w_2^h,\\\cM^-_{\fL_0} w_1^h\le \cM^-_{\fL^1_0}
w_1^h&=\cM^-_{\fL^1_0} [w^h-w_2^h]\le 0-\cM^-_{\fL^1_0}
w_2^h.\end{split}
\end{equation}
In order to apply Theorem \ref{thm-9.2}, we have only to show that
$|\cM^+_{\fL^1_0} w_2^h|$ and $|\cM^-_{\fL^1_0} w_2^h|$ are bounded
on $B_{r-\dt/2}$ by $C\,\|u\|_{L^{\iy}(\BR^n)}$ for some universal
constant. To show this, we prove that it is true for any operator
$\cL^t\in\fL^1_0$. Take any $\cL^t\in\fL^1_0$. Since $w_2^h\equiv 0$
on $B_{r-\dt/2}$, the jump part of $\cL^t w_2^h$ vanishes on
$B_{r-\dt/2}$. Thus by the change of variable and the mean value
theorem we have that
\begin{equation*}\begin{split}\cL^t w_2^h(x)
&=\int_{\BR^n}\f{\bigl(1-\tau_h\phi(x+y)\bigr)\tau_h u(x+y)-\bigl(1-\phi(x+y)\bigr)\,u(x+y)}{|h|^{\bt}}\,K(y)\,dy\\
&\qquad+\int_{\BR^n}\f{\bigl(\tau_h\phi(x+y)-\phi(x+y)\bigr)\tau_h u(x+y)}{|h|^{\bt}}\,K(y)\,dy\\
&=\int_{\BR^n}\bigl(1-\phi(x+y)\bigr)\,u(x+y)\,\f{K(y-h)-K(y)}{|h|^{\bt}}\,dy\\
&\qquad+\int_{\BR^n}\f{\bigl(\int_0^1[\n\phi(x+y+\tau h)\cdot
h]\,d\tau\bigr)\tau_h u(x+y)}{|h|^{\bt}}\,K(y)\,dy\\&\fd\cL^t_1
w_2^h(x)+\cL^t_2 w_2^h(x)\end{split}\end{equation*} for any $x\in
B_{r-\dt/2}$. Set $\vr_1=\dt/8$. Since $(1-\phi(x+y))u(x+y)=0$ for
any $x\in B_{r-\dt/2}$ and $|y|<\dt/8$, by \eqref{eq-10.1} we obtain
that
\begin{equation}\label{eq-10.3}
\begin{split}\bigl|\cL^t_1
w_2^h(x)\bigr|&\le|h|^{1-\bt}\,\|u\|_{L^{\iy}(\BR^n)}\int_{y\in
B^c_{\vr_1}}\f{|K(y)-K(y-h)|}{|h|}\,dy\le C\,\|u\|_{L^{\iy}(\BR^n)}
\end{split}\end{equation} for any $x\in B_{r-\dt/2}$ and $|h|<\dt/16$. Since
$\tau_h\phi(x+y)=\phi(x+y)$ for every $x\in B_{r-\dt/2}$,
$|y|<\dt/8$ and $|h|<\dt/16$, we have that
\begin{equation}\label{eq-10.4}\bigl|\cL^t_2
w_2^h(x)\bigr|\le|h|^{1-\bt}\|\n\phi\|_{L^{\iy}(\BR^n)}\,\|u\|_{L^{\iy}(\BR^n)}\int_{\BR^n\s
B_{\vr_1}}\f{1}{|y|^{n+\sm}}\,dy\le
C\,\|u\|_{L^{\iy}(\BR^n)}
\end{equation} for any $x\in B_{r-\dt/2}$ and
$|h|<\dt/16$. Then it follows from \eqref{eq-10.2}, \eqref{eq-10.3}
and \eqref{eq-10.4} that $\cM^+_{\fL_0^1}
w_1^h\ge-C\,\|u\|_{L^{\iy}(\BR^n)}$ and $\cM^-_{\fL_0^1} w_1^h\le
C\,\|u\|_{L^{\iy}(\BR^n)}$ on $B_{r-\dt/2}$ for any $h\in
B_{\dt/16}$. Thus by Theorem \ref{thm-9.2} we have that
\begin{equation}\label{eq-10.5}\begin{split}\|w^h\|_{\rC^{\ap}(B_{r-\dt})}&=\|w_1^h\|_{\rC^{\ap}(B_{r-\dt})}
\le C\,\|w_1^h\|_{L^{\iy}(B_{r-\dt/2})}+C\,\|u\|_{L^{\iy}(\BR^n)}\\
&\le
C\,\|u\|_{\rC^{0,\bt}(B_{r-\dt/4})}+C\,\|u\|_{L^{\iy}(\BR^n)}\end{split}\end{equation}
for any $h\in B_{\dt/16}$. By the standard telescopic sum argument
\cite{CC}, we obtain that $\|u\|_{\rC^{\ap+\bt}(B_{r-\dt})}\le
C\,(\,\|u\|_{L^{\iy}(\BR^n)}+|\cI^{\pm} 0|).$ Repeating the above
argument, after $[1/\ap]$-steps we have that
\begin{equation}\label{eq-10.6}\|u\|_{\rC^{0,1}(B_{3/4})}\le
C\,(\,\|u\|_{L^{\iy}(\BR^n)}+|\cI^{\pm} 0|).\end{equation} For any
unit vector $e\in S^{n-1}$, we consider the following incremental
quotients with the same reasoning
$$w^t_e(x)=\f{\tau_{te}u(x)-u(x)}{h},\,t>0.$$ If we choose $r=5/6$, $\dt=1/3$
and $\bt=1$ in \eqref{eq-10.5}, then by \eqref{eq-10.5} and
\eqref{eq-10.6} we obtain that $\|w^t_e\|_{\rC^{\ap}(B_{1/2})}\le
C\,(\,\|u\|_{L^{\iy}(\BR^n)}+|\cI^{\pm} 0|)$ for any unit vector
$e\in S^{n-1}$ and for any $t$ with $|t|<1/48$. From this, taking
$t\downarrow 0$ we conclude that $\|u\|_{\rC^{1,\ap}(B_{1/2})}\le
C\,(\,\|u\|_{L^{\iy}(\BR^n)}+|\cI^{\pm} 0|).$ Hence we complete the
proof. \qed

\noindent{\bf Acknowledgement.} This work had been started during
Yong-Cheol Kim was visiting to University of Texas at Austin in the fall
semester 2008 for his sabbatical year. He would like to thank for
kind hospitality and deep concern of Professor Luis A. Caffarelli
during that time. Ki-Ahm Lee was supported by the Korea Research
Foundation Grant funded by the Korean Government(MOEHRD,
Basic Research Promotion Fund)( KRF-2008-314-C00023).

\end{document}